\input ./style/arxiv-general.cfg
\documentclass[aop,MSNbibl,seceqn,dvips]{arximspdf}
\makeatletter
   \@ifpackageloaded{graphicx}{}{\usepackage{graphicx}}
\makeatother
\usepackage{mathrsfs}

\doi{10.1214/15-AOP1020}
\volume{44}
\issue{3}
\pubyear{2016}
\firstpage{2147}
\lastpage{2197}
\docsubty{FLA}

\makeatletter
\newcommand{\rrvert}{\vert}
\newcommand{\rrVert}{\Vert}
\newcommand{\llvert}{\vert}
\newcommand{\llVert}{\Vert}
\renewcommand{\mid}{|}
\newcommand{\eqref}[1]{(\ref{#1})}
\newtheorem{theorem}{Theorem}[section]
\newtheorem{proposition}[theorem]{Proposition}
\newtheorem{lemma}[theorem]{Lemma}
\newtheorem{corollary}[theorem]{Corollary}
\newproclaim{remark}[theorem]{Remark}
\newproclaim{example}[theorem]{Example}
\newcommand{\dist}{\operatorname{dist}}
\newcommand{\dom}{\operatorname{dom}}
\def\AA{\mathbb{A}}
\def\GG{\mathbb{G}}
\def\NN{\mathbb{N}}
\def\QQ{\mathbb{Q}}
\def\RR{\mathbb{R}}
\def\SS{\mathbb{S}}
\def\XX{\mathbb{X}}
\def\YY{\mathbb{Y}}
\def\ZZ{\mathbb{Z}}
\def\a{\alpha}
\def\b{\beta}
\def\g{\gamma}
\def\d{\delta}
\def\k{\kappa}
\def\l{\lambda}
\def\s{\sigma}
\def\o{\omega}
\def\bC{\mathbf{C}}
\def\bE{\mathbf{E}}
\def\bK{\mathbf{K}}
\def\bL{\mathbf{L}}
\def\bM{\mathbf{M}}
\def\bP{\mathbf{P}}
\def\bQ{\mathbf{Q}}
\def\cB{\mathcal{B}}
\def\cC{\mathscr{L}_1}
\def\cH{\mathcal{H}}
\def\cN{\mathcal{N}}
\def\cY{\mathcal{Y}}
\def\cX{\mathcal{X}}
\def\N{\mathrm{N}}
\def\sC{\mathscr{C}}
\def\dint{\mathrm{d}}
\def\dTV{\mathrm{d}_{\mathrm{TV}}}
\def\dTVRandom{\mathbf{d}_{\mathbf{TV}}}
\def\dKR{\mathbf{d}_{\mathbf{KR}}}
\def\dW{\mathbf{d}_{\mathbf{W}}}
\def\dK{\mathbf{d}_{\mathbf{K}}}
\def\vol{\operatorname{vol}}
\def\dist{\operatorname{dist}}
\makeatother

\begin{document}
\begin{frontmatter}

\title{Functional Poisson approximation in Kantorovich--Rubinstein
distance with applications to U-statistics and stochastic geometry}
\runtitle{Functional Poisson approximation}

\begin{aug}
\author[A]{\fnms{Laurent}~\snm{Decreusefond}\thanksref{T1}\ead[label=e1]{laurent.decreusefond@telecom-paristech.fr}},
\author[B]{\fnms{Matthias}~\snm{Schulte}\thanksref{T2}\ead[label=e2]{matthias.schulte@kit.edu}}
\and
\author[C]{\fnms{Christoph}~\snm{Th\"ale}\corref{}\thanksref{T3}\ead[label=e3]{christoph.thaele@rub.de}}
\runauthor{L. Decreusefond, M. Schulte and C. Th\"ale}
\affiliation{Telecom ParisTech, Karlsruhe Institute of
Technology and Ruhr University Bochum}
\address[A]{L. Decreusefond\\
Telecom ParisTech\\
Rue Barrault 46\\
F-75634 Paris cedex 13\\
France\\
\printead{e1}}
\address[B]{M. Schulte\\
Department of Mathematics\\
Institute of Stochastics\\
Karlsruhe Institute of Technology\\
D-76128 Karlsruhe\\
Germany\\
\printead{e2}}
\address[C]{C. Th\"ale\\
Faculty of Mathematics\\
Ruhr University Bochum\\
NA 3/68, D-44780 Bochum\\
Germany\\
\printead{e3}}
\end{aug}
\thankstext{T1}{Supported in part by ANR Masterie.}
\thankstext{T2}{Supported by the German Research Foundation (DFG) through the research unit ``Geometry and Physics of Spatial Random Systems'' Grant HU 1874/3-1.}
\thankstext{T3}{Supported by the German research foundation (DFG) via SFB-TR 12.}

%
\received{\smonth{6} \syear{2014}}
%
\revised{\smonth{3} \syear{2015}}

%
\begin{abstract}
A Poisson or a binomial process on an abstract state space and a
symmetric function $f$ acting on
$k$-tuples of its points are considered. They induce a point process on
the target space of $f$. The main result is a functional limit theorem
which provides an upper bound for an optimal transportation distance
between the image process and a Poisson process on the target space.
The technical background are a version of Stein's method for Poisson
process approximation, a Glauber dynamics representation for the
Poisson process and the Malliavin formalism. As applications of the
main result, error bounds for approximations of U-statistics by
Poisson, compound Poisson and stable random variables are derived, and
examples from stochastic geometry are investigated.
\end{abstract}

%
\end{frontmatter}

\section{Introduction}\label{sec1}

The arguably most prominent functional limit theorem is Donsker's
invariance principle. It asserts that the distribution of a linear
interpolation between the points of a suitably re-scaled random walk
converges to the Wiener measure on the space of continuous functions
on $\RR_+$, the nonnegative real half-line; see, for example, \cite
{KallenbergFoundations}, Corollary 16.7. Besides the Wiener process,
there is
another fundamental stochastic process, which plays an important
role in many branches of probability theory and its applications,
namely the Poisson
process. However, functional limit theorems involving the Poisson
process have found much less attention in the literature. The aim of
this paper is to provide a quantitative version of a functional limit
theorem for Poisson processes
and to derive from it error bounds for the probabilistic approximation
of U-statistics by a Poisson, a compound Poisson
or a stable random variable. We demonstrate the versatility of our
results by applying these bounds to
functionals of random geometric graphs, distance-power statistics,
nonintersecting flat processes and random polytopes.

Let us informally describe the set-up of this paper; precise
definitions and statements follow in Section~\ref{secmain-result}.
Let $(\XX,\cX)$ and $(\YY,\cY)$ be two measurable spaces
(\mbox{satisfying} some mild regularity assumptions, see below), let $\bK_1$
be a probability measure on $\XX$ and fix an integer $k\geq1$.
Moreover, for each $n\in\NN$ let $f_n\dvtx\dom f_n \to\YY$ be a
symmetric mapping whose domain $\dom f_n$ is a symmetric subset of $\XX
^k$. Next, consider a collection $\beta_n=\{X_1,\ldots,X_n\}$ of
$n\geq k$ i.i.d. random elements $X_1,\ldots,X_n$ of $\XX$ with
distribution $\bK_1$. We apply for each $n\geq k$, $f_n$ to every
$k$-tuple of distinct elements of $\beta_n$. This induces a point
process $\xi_n$ on $\YY$ of the form
\[
\xi_n=\frac{1}{k!}\sum_{(x_1,\ldots,x_k)\in\beta^k_{n,\neq}\cap
\dom f_n}
\delta_{f_n(x_1,\ldots,x_k)}, %
\]
where $\beta_{n,\neq}^k=\{(x_1,\ldots,x_k)\in\beta_n^k\dvtx x_j\neq
x_j, i\neq j, i,j=1,\ldots,k\}$ and $\delta_y$ stands for the unit
mass Dirac measure concentrated at $y\in\YY$.

The motivation for studying the point processes $\xi_n$ as defined
above comes from the theory of U-statistics and from a class of extreme
value problems arising in stochastic geometry. At first, if $\dom
f_n=\XX^k$ and $\YY=\RR$, the points of $\xi_n$ can be regarded as
the summands of the U-statistic
\[
S_n=\frac{1}{k!} \sum_{(x_1,\ldots,x_k)\in\beta_{n,\neq}^k}
f_n(x_1,\ldots,x_k) . %
\]
These objects play a prominent role in large parts of probability
theory and mathematical statistics, and an analysis of the point
process of summands is helpful for the understanding of their
(asymptotic) properties. On the other hand, in several problems arising
in stochastic geometry, one is interested in extreme values of the type
\[
\min_{(x_1,\ldots,x_k)\in\beta_{n,\neq}^k} f_n(x_1,
\ldots,x_k) %
\]
in case that $\dom f_n=\XX^k$ and $\YY=[0,\infty)$. Clearly, this
minimum is the distance from the origin to the first point of the point
process $\xi_n$. For these reasons, a study of the point processes
$\xi_n$ unifies both mentioned problems.

The intensity measure $\bL_n$ of $\xi_n$ is given by
%
\begin{eqnarray}
\bL_n(A)=\bE\xi_n(A)=\frac{(n)_k}{k!}\int
_{\dom f_n }{\mathbf1}\bigl(f_n(x_1,
\ldots,x_k)\in A\bigr) \bK_1^k\bigl(\dint(x_1,\ldots,x_k)\bigr),\nonumber
\\
\eqntext{A\in\cY,}
\end{eqnarray}
where $(n)_k$ is the descending factorial. Our main result, Theorem
\ref{thmMain} below, provides an upper bound for the
Kantorovich--Rubinstein distance $\dKR(\xi_n,\zeta)$ between $\xi
_n$ and a Poisson process $\zeta$ on $\YY$ with finite intensity
measure $\bM$. Here, the Kantorovich--Rubinstein distance is a variant
of an optimal transportation distance, which measures the closeness
between two point processes or, more precisely, their distributions. In
particular, we show that $\xi_n$ converges in Kantorovich--Rubinstein
distance, and thus in distribution to $\zeta$ if
\[
\dTV(\bL_n,\bM)\to0\quad\mbox{and} \quad\bE\xi_n(
\YY)^2\to\bM(\YY)^2+\bM(\YY)\qquad\mbox{as } n\to\infty,
\]
where $\dTV( \cdot, \cdot)$ denotes the total variation
distance of measures on $\YY$. More precisely, the upper bound for the
Kantorovich--Rubinstein distance only depends on $\dTV(\bL_n,\bM)$
and the first two moments of $\xi_n(\YY)$. This is a functional
version of the famous results by Arratia, Goldstein and Gordon \cite{ArratiaGoldsteinGordon}, and Chen \cite{Chen75} that ``two moments suffice for Poisson
approximation.''

Besides the binomial process $\beta_n$ of $n$ independent and
identically distributed points, we also allow the input process to be a
Poisson process on $\XX$ with a \mbox{$\sigma$-}finite intensity measure. 
In some instances, an underlying Poisson process is more natural and
sometimes even unavoidable, especially if the underlying point process
on $\XX$ is supposed to have infinitely many points. To exploit this
flexibility, we consider both set-ups in parallel.

Poisson process approximation has been studied by several authors by
means of Stein's method, but to the best of our knowledge this is the
first paper where the Kantorovich--Rubinstein distance is investigated.
The works of Barbour \cite{Barbour88}, Barbour and Brown \cite
{BarbourBrown} and the last chapter of the monograph \cite
{BarbourHolstJanson} of Barbour, Holst and Janson concern Poisson
process approximation in the total variation distance. But since the
total variation distance is not suitable for all problems and since the
so-called \textit{Stein magic factors} do not get small if $\bL_n(\YY)$
is large (in contrast to classical Poisson approximation), one often
uses weaker notions of distance. Starting with the work of Barbour and
Brown \cite{BarbourBrown} and Barbour, Holst and Janson \cite
{BarbourHolstJanson}, this has been done by Brown, Chen, Schuhmacher,
Weinberg and Xia \cite
{BrownWeinbergXia,BrownXia1995,BrownXia2001,ChenXia,Schuhmacher2005,Schuhmacher2009,SchuhmacherXia}.
Our work goes in the opposite direction since the
Kantorovich--Rubinstein distance between point processes is stronger
than the total variation distance in the sense that convergence in
Kantorovich--Rubinstein distance implies convergence in total variation
distance, but not vice versa. Roughly speaking and in a transferred
sense, the Kantorovich--Rubinstein distance is related to the total
variation distance between point processes as the Wasserstein distance
is related to the total variation distance for integer-valued random
variables. Since its test functions are allowed to take values other
than zero and one, the Kantorovich--Rubinstein distance is more
sensitive to the behavior and the geometry of the compared point
processes than the total variation distance. Let us further remark that
in the recent paper \cite{SchuhmacherStucki}, Schuhmacher and Stucki
consider the total variation distance between two Gibbs processes. This
includes Poisson process approximation as a special case. However, the
approximated point processes of the present paper do not, in general,
satisfy the technical conditions assumed in \cite{SchuhmacherStucki}
since they are not necessarily hereditary.

Besides the notion of distance and its connection to the theory of
optimal transportation, the other main ingredient of our approach is a
functional version of Stein's method for Poisson process approximation.
It relies on a Glauber dynamics representation for Poisson processes
and the Malliavin formalism. More precisely, we use an
integration-by-parts argument on the target space and then a
commutation relation between the discrete gradient on that space and
the semi-group associated with the Glauber dynamics. This way we avoid
the explicit computation and investigation of a solution of the Stein equation.
We would like to highlight that our approach is generic and depends
only on the underlying random structure (here, a binomial or a Poisson
process) and not on a very specific model so that extensions to other
probabilistic frameworks (such as Gaussian random measures or
Rademacher sequences) should also be possible. However, they are beyond
the scope of this paper and will be treated elsewhere.

To demonstrate the versatility of our new functional limit theorem, we
consider probabilistic approximations of U-statistics over binomial or
Poisson input processes. In a first regime, we consider the Poisson
approximation of U-statistics and provide an error bound for the
Wasserstein distance. Our result improves and extends earlier works of
Barbour and Eagleson \cite{BarbourEagleson} and Peccati \cite
{Peccati11ChenStein}. The second regime concerns compound Poisson
approximation of U-statistics in total variation distance. Here, we do
not impose any conditions on the nature of the compound Poisson
distribution, which is allowed to be discrete or continuous. In
contrast, previous results for the compound Poisson approximation via
Stein's method only deal with the discrete case; see, for example, the
work of Barbour, Chen and Loh \cite{BarboudChenLoh}, the survey \cite
{BarbourChryssaphinou} of Barbour and Chryssaphinou and especially the
paper \cite{EichelsbacherRoos} of Eichelsbacher and Roos, who consider
U-statistics over a binomial input process. In this light, we
generalize the results of \cite{EichelsbacherRoos} to a larger class
of limiting distributions and also to the case of an underlying Poisson
process. In a third regime, we use our functional limit theorem to
investigate probabilistic approximations of U-statistics by $\alpha
$-stable random variables with $0<\alpha<1$ and to derive explicit
error bounds for the Kolmogorov distance. In their previous work \cite
{DehlingEtAl}, Dabrowski et al. also obtained $\alpha$-stable limits
for U-statistics from point process convergence results. However, their
technique does not allow any conclusion about a rate of convergence.

Finally, we apply our general result to problems arising in stochastic
geometry. Random geometric graphs are one of the fundamental models of
spatial stochastics; see \cite{PenroseBook}, for example. We derive
limit theorems for several U-statistics of random geometric graphs,
where the limiting distributions are Poisson or compound Poisson, and
show a new point process limit theorem for the midpoints of short
edges. As further examples, we consider distance-power statistics with
$\alpha$-stable limit distributions, midpoints between
nonintersecting Poisson $m$-flats which are close together and the
diameter of random polytopes with vertices on the sphere.

In a natural way our paper continues the line of research on point
process convergence and extreme values initiated by the second and the
third author in \mbox{\cite{STScalingLimits,STFlats}}, where the proofs are
based on the main result of \cite{Peccati11ChenStein} and the
underlying point process has to be Poisson. In contrast to these
previous works our technique also allows us to deal with an underlying
binomial process and delivers in both cases bounds for the
Kantorovich--Rubinstein distance. Furthermore, the bounds derived here
improve the rates of convergence of some of the scalar limit theorems
from \cite{STScalingLimits,STFlats}. Our findings also complement the
works \cite{DecreusefondWassersteinDistance} and \cite
{DecreusefondRubinstein} of the first author with Joulin and Savy,
concerning the Kantorovich--Rubinstein distance on configuration spaces
and related notions.

This paper is organized as follows. Before we present our main result
for Poisson process convergence in Section~\ref{secmain-result}, we
recall in Section~\ref{secPreliminaries} some necessary notation and
results about point processes and also summarize some facts from convex
geometry which are important for our examples from stochastic geometry.
The proof of our main result in Section~\ref{secProof-main-result} is
prepared by a brief discussion of the underlying Stein principle in
Section~\ref{secgener-stein-princ} and the Glauber dynamics, a key
step in our argument, in Section~\ref{secglaub-dynam-poiss}.
Section~\ref{secApplications} is devoted to applications of our
functional limit theorem to probabilistic approximations of
U-statistics and to problems from stochastic geometry.

\section{Preliminaries}\label{secPreliminaries}

In the present section we introduce some basic notions and notation,
which are used in the text. Throughout, $(\Omega,{\mathcal F},\bP)$
will be
an abstract probability space, which is rich enough to carry all the
random objects we deal with. Expectation with respect to $\bP$ is
denoted by $\bE$.

\subsection{Configuration spaces}
Let $(\YY,\cY)$ be a lcscH space; that is, $\YY$ is a
topological space with countable base such that every point in $\YY$
has a compact neighborhood and such that any two points of $\YY$ can
be separated by disjoint neighborhoods. Such a space is separable and
completely metrizable. Here, $\cY$ denotes the Borel $\sigma$-field
generated by the topology of $\YY$. By $\N_\YY$ we denote the space
of $\s$-finite
counting measures (i.e., point configurations) on $\YY$, whereas
$\widetilde{\N}_\YY$ and $\widehat{\N}_\YY$ stand for the sets of
all locally finite (i.e., bounded on all relatively compact sets) and
finite counting measures on $\YY$, respectively. By a slight abuse
of notation we will write $y\in\o$ if $y\in\YY$ is charged by the
measure $\o$ and also use the set-notation $\o_1\subset\o_2$ to
indicate that $\o_1$ is a sub-configuration of $\o_2$ (with a similar
meaning we also understand $\o_2\setminus\o_1$). Let $\cN_\YY$ be
the $\s$-field on $\N_\YY$ generated by
the mappings
\[
\psi_A\dvtx\N_\YY\to\NN_0\cup\{\infty\} ,\qquad \o
\mapsto\o(A) , A\in\cY,
\]
where $\NN_0:=\NN\cup\{0\}$ is the set of natural numbers including
zero. We equip $\widetilde{\N}_\YY$ and $\widehat{\N}_\YY$ with
the corresponding trace $\sigma$-fields of $\cN_\YY$. The $\sigma
$-field of $\widetilde{\N}_\YY$ is then the Borel $\s$-field for the
vague topology on $\widetilde{\N}_\YY$, which is generated by the
mappings
\[
e_g\dvtx\widetilde{\N}_\YY\to[0,\infty),\qquad \o\mapsto\int
_{\YY} g \,\dint\o,
\]
where $g\geq0$ is a continuous function on $\YY$ with compact
support, and the space $\widetilde{\N}_\YY$ equipped with the vague
topology becomes a Polish space; see Theorem A2.3 in \cite
{KallenbergFoundations}. A point process (or random counting measure)
$\mu$ is a random element in~$\N_\YY$. By a locally finite point
process and a finite point process, we mean random elements in
$\widetilde{\N}_\YY$ and $\widehat{\N}_\YY$, respectively. It
follows from \cite{SW}, Lemma 3.1.3, that a point process $\mu$ can
almost surely be represented as
\[
\mu=\sum_{i=1}^{\mu(\YY)} \delta_{x_i}
\qquad\mbox{with } x_i\in\YY,  i\in\NN\mbox{ and }\mu(
\YY)\in\NN_0\cup\{ \infty\} , %
\]
where $\d_y$ stands for the unit mass Dirac measure concentrated at
$y\in\YY$. Thus we may interpret $\mu$ also as a random collection
of points, taking into account potential multiplicities.

\subsection{Poisson processes}
Let $\bM$ be a $\s$-finite measure
on $\YY$, and let $\bM^k$ stand for its $k$-fold product measure. By
a Poisson process on $\YY$ with intensity measure~$\bM$, we
understand a point process $\zeta$ with the properties that:
(i) for any $B\in\cY$, the random variable $\zeta(B)$ is Poisson
distributed with mean $\bM(B)$ and (ii) $\zeta$ is independently
scattered; that is, for any $n\in\NN$ and disjoint $B_1,\ldots
,B_n\in\cY$
the random variables $\zeta(B_1),\ldots,\zeta(B_n)$ are
independent. We notice that if $\bM$ is a finite measure, $\zeta$
charges almost surely only a finite number of points in $\YY$, whose
total number follows a Poisson distribution with mean $\bM(\YY)$. We
will write $\bP_\zeta$ for the distribution of $\zeta$ on $\N_\YY
$. In
this paper we will speak about a homogeneous Poisson process on a set
$A\in\mathcal{B}(\RR^d)$, where $\mathcal{B}(\RR^d)$ is the Borel
$\sigma$-field on $\RR^d$,
if the intensity measure is a multiple of the restriction of the
Lebesgue measure to $A$. Also, if $d=1$, a homogeneous Poisson
process $\zeta$ on $[0,\infty)$ can be thought of as a piecewise
deterministic (pure jump) stochastic process in continuous time,
starting at zero and having jumps of size one and i.i.d. exponentially
distributed waiting times between the jumps. The points of
discontinuity of this random process are the jump times of $\zeta$.

One of our main tools to deal with Poisson functionals (by these we
mean real-valued random variables depending only on a Poisson process)
is the
multivariate Mecke formula \cite{SW}, Corollary 3.2.3, which says that
for any
integer $k\geq1$ and any measurable and nonnegative $f\dvtx\YY^k\times
\N_\YY\to\RR$,
%
\begin{eqnarray}\label{eqMecke}
&& \bE\sum_{(y_1,\ldots,y_k)\in\zeta_{\neq
}^k}f(y_1,\ldots
,y_k,\zeta)
\nonumber\\[-8pt]\\[-8pt]\nonumber
&&\qquad =\int_{\YY^k}\bE f(y_1,
\ldots,y_k,\zeta+\d_{y_1}+\cdots+\d_{y_k})
\bM^k \bigl(\dint(y_1,\ldots,y_k) \bigr) ,
\end{eqnarray}
where $\zeta_{\neq}^k$ is the collection of all $k$-tuples of distinct
points charged by $\zeta$. If the point process $\zeta$ is simple
[i.e., if $\zeta(\{y\})\in\{0,1\}$ almost surely for any $y\in\YY
$], $\zeta_{\neq}^k$~can be written as
\[
\zeta_{\neq}^k=\bigl\{(y_1,
\ldots,y_k)\in\YY^k\dvtx y_i\neq y_j
\in\zeta\mbox{ for } i\neq j, i,j=1,\ldots,k\bigr\} , %
\]
while in the nonsimple case distinct points can have the same location.
We remark that \eqref{eqMecke} with $k=1$ is even a characterizing
property of the Poisson process $\zeta$; cf. Theorem 3.2.5 of \cite{SW}.

\subsection{Binomial processes}
Let $\bM_1$ be a probability measure on $\YY$. A binomial process with
intensity measure $\bM:=n\bM_1$, $n\in\N$, is a collection of $n$
random points,
distributed independently according to the measure $\bM_1$. This
process also arises by conditioning a Poisson process with intensity
measure $\bM$ on having exactly $n$ points. In this paper we shall
denote the random counting measure induced by such a binomial process
by $\beta_n$. We also write $\beta_{n,\neq}^k$ to indicate
the collection of all $k$-tuples of distinct points charged by
$\beta_n$. Then the counterpart to the multivariate Mecke formula
\eqref{eqMecke} for a
binomial process reads as follows:
%
\begin{eqnarray}
\label{eqMomentMeasureBinomial}
&& \bE\sum_{(y_1,\ldots,y_k)\in\beta
^k_{n,\neq}} f(y_1,
\ldots,y_k,\beta_n)
\nonumber\\[-8pt]\\[-8pt]\nonumber
&&\qquad = (n)_k\int_{\YY^k}\bE f(y_1,
\ldots,y_k,\beta_{n-k}+\d_{y_1}+\cdots+
\d_{y_k}) \bM_1^k \bigl(\dint(y_1,
\ldots,y_k) \bigr) ,\hspace*{-25pt}
\end{eqnarray}
where $(n)_k:=n(n-1)\cdots(n-k+1)$ is the descending factorial and
$f$ is a real-valued nonnegative measurable function on $\YY^k\times
\N_\YY$. This can easily be
seen directly and is also a special case of the Georgii--Nguyen--Zessin
formula, for which we refer to \cite{DaleyVereJones2008}, Proposition 15.5.II.

\subsection{Probability distances}
In order to compare two real-valued random variables $Y_1$ and $Y_2$
(or more precisely their distributions) and to measure their closeness,
we use several probability distances in this paper. The Kolmogorov
distance of $Y_1$ and $Y_2$ is given by
\[
\dK(Y_1,Y_2):=\sup_{z\in\RR}\bigl\llvert
\bP(Y_1\leq z)-\bP(Y_2\leq z)\bigr\rrvert, %
\]
while the total variation distance is
\[
\dTVRandom(Y_1,Y_2):=\sup_{A\in\cB(\RR)}\bigl
\llvert\bP(Y_1\in A)-\bP(Y_2\in A)\bigr\rrvert,
\]
where, recall, $\cB(\RR)$ stands for the Borel $\sigma$-field on
$\RR$. If $Y_1$ and $Y_2$ are integer-valued random variables, we can
re-write their total variation distance as
\[
\dTVRandom(Y_1,Y_2) = \frac{1}{2} \sum
_{k\in\ZZ} \bigl\llvert\bP(Y_1=k)-\bP
(Y_2=k)\bigr\rrvert. %
\]
Let us denote by $\operatorname{Lip}(1)$ the set of all functions
$h\dvtx\RR\to\RR
$ whose Lipschitz constant is at most one and define the Wasserstein
distance of two real-valued random variables $Y_1$ and $Y_2$ by
\[
\dW(Y_1,Y_2):=\sup_{h\in\operatorname{Lip}(1)}\bigl\llvert
\bE h(Y_1)-\bE h(Y_2) \bigr\rrvert. %
\]
These probability distances all have the property that they imply
convergence in distribution, meaning that for a sequence $(Y_n)_{n\in
\NN}$ of random variables convergence in distribution to another
random variable $Y$ is implied by
%
\begin{equation}
\label{eqConvergenceDistances} \mathbf{d_I}
(Y_n,Y)\to0
\qquad
\mbox{as }n\to\infty,\mbox{ for some }\mathbf{I}\in\{\mathbf
{K},\mathbf{TV},\mathbf{W}\} .
\end{equation}
Moreover, for integer-valued random variables $Y_1$ and $Y_2$, let us
mention the general inequality
%
\begin{equation}
\label{eqInequalityDistances} \dK(Y_1,Y_2)\leq\dTVRandom(Y_1,Y_2)
\leq\dW(Y_1,Y_2) ,
\end{equation}
which directly follows from the definitions of the involved probability
distances and the fact that $Y_1$ and $Y_2$ are concentrated on the
integers. Note that \eqref{eqInequalityDistances} does not remain
valid for general real-valued random variables.

\subsection{Kantorovich--Rubinstein distance}\label{sec2.5}
We define the total variation distance between two measures $\nu_1$
and $\nu_2$ on $\YY$
by
\[
\dTV(\nu_1,\nu_2):=\mathop{\sup_{A\in\cY}}_{\nu
_1(A),\nu_2(A)<\infty
}\bigl
\llvert\nu_1(A)-\nu_2(A)\bigr\rrvert, %
\]
a notion that should not be confused with the total variation distance
between random variables introduced above. Note that $\dTV(\nu_1,\nu
_2)$ can in principle take any value in $[0,\infty]$.

We say that a map $h\dvtx\widetilde{\N}_\YY\to\RR$ is $1$-Lipschitz if
\[
\bigl\llvert h(\o_1)-h(\o_2)\bigr\rrvert\leq\dTV(
\o_1,\o_2)\qquad\mbox{for all } \o_1,
\o_2\in\widetilde{\N}_\YY, %
\]
and denote by $\cC$ the set of all these maps which are measurable.

The Kantorovich--Rubinstein distance between two probability measures
$\bQ_1$ and $\bQ_2$ on $\N_\YY$
is defined as the optimal transportation cost
%
\begin{equation}
\label{eqRubinsteinDistance} \dKR(\bQ_1,\bQ_2):=\inf
_{\bC\in\Sigma(\bQ_1,\bQ_2)}\int_{\N
_\YY\times\N_\YY}\dTV(\o_1,
\o_2) \bC\bigl(\dint(\o_1,\o_2)\bigr)
\end{equation}
for the cost function $\dTV( \cdot, \cdot)$, where $\Sigma
(\bQ_1,\bQ_2)$ denotes the set of probability measures on
$\N_\YY\times\N_\YY$ with first marginal $\bQ_1$ and the second marginal
$\bQ_2$ (i.e., couplings of $\bQ_1$ and $\bQ_2$). If $\bQ_1$ and
$\bQ_2$ are concentrated on $\widetilde{\N}_\YY$, there is at least
one coupling $\bC\in\Sigma(\bQ_1,\bQ_2)$ for which the infimum in
\eqref{eqRubinsteinDistance} is attained according to \cite
{Villani2009}, Theorem~4.1, and the Kantorovich duality theorem
\cite{Villani2009}, Theorem~5.10, says that this minimum equals
%
\begin{equation}
\label{eqRubinsteinDistanceEquivalent} \dKR(\bQ_1,\bQ_2)=\sup\biggl
\llvert\int
_{\widetilde{\N}_\YY} h(\omega) \bQ_1(\dint\omega) -\int
_{\widetilde{\N}_\YY} h(\omega) \bQ_2(\dint\omega) \biggr\rrvert
,
\end{equation}
where the supremum is over all $h\in\cC$ that are integrable with
respect to $\bQ_1$~and~$\bQ_2$. 

By abuse of notation we will also write
$\dKR(\zeta_n,\zeta)$ instead of $\dKR(\bQ_n,\bQ)$ if the point process
$\zeta_n$ on $\YY$ has distribution $\bQ_n$ for any $n\geq1$ and
the point process
$\zeta$ on $\YY$ has distribution $\bQ$. Note that the integrability
condition in \eqref{eqRubinsteinDistanceEquivalent} is automatically
fulfilled for all $h\in\cC$ if $\bE\zeta_n(\YY)<\infty$ and $\bE
\zeta(\YY)<\infty$. The Kantorovich--Rubinstein distance is also
called Wasserstein distance, Monge--Kantorovich \mbox{distance} or Rubinstein
distance. For a detailed discussion of the terminology we refer to the
bibliographic notes of Chapter~6 in \cite{Villani2009}.

The following result ensures that convergence of locally finite point
processes in Kantorovich--Rubinstein distance implies convergence in
distribution.

%
\begin{proposition}\label{thmPPPconv1}
Assume that $(\zeta_n)_{n\in\NN}$ is a sequence of locally finite
point processes on $\YY$ and that $\zeta$ is another locally finite
point process on
$\YY$ such that $\dKR(\zeta_n,\zeta)\to0$, as $n\to\infty$. Then
$\zeta_n$ converges in distribution to $\zeta$, as $n\to\infty$.
\end{proposition}

\begin{pf}
The structure of the vague topology on $\widetilde{\N}_\YY$ implies
that it is
necessary and sufficient to prove that for any continuous $g \dvtx
\YY\to\RR$ with compact support, the random variables $\int g \,\dint
\zeta_n$ converge in distribution to $\int g \,\dint\zeta$; see
\cite{KallenbergFoundations}, Theorem~16.16.
By \eqref{eqConvergenceDistances}, it is sufficient to
show that for all Borel sets $B\subset\RR$, we have that
\[
\bE e_{g,B}(\zeta_n)\to\bE e_{g,B}(\zeta) \qquad
\mbox{as } n\to\infty,
\]
where $e_{g,B}\dvtx\widetilde{\N}_\YY\to\RR, \omega\mapsto{\mathbf
1}(\int g \,\dint\omega\in B)$.
To show this, we notice that for each $g$ and $B$ as above the mapping
$e_{g,B}$ belong to $\cC$, whence
\[
\bigl\llvert\bE e_{g,B}(\zeta_n)-\bE e_{g,B}(
\zeta)\bigr\rrvert\leq\dKR(\zeta_n,\zeta),
\]
and the result follows.
\end{pf}

An alternative distance to measure the closeness of two point processes
$\zeta_1$ and $\zeta_2$ on $\YY$ is the total variation distance
\[
\dTVRandom(\zeta_1,\zeta_2):=\sup_{A\in\cN_\YY}
\bigl\llvert\bP(\zeta_1\in A)-\bP(\zeta_2\in A)\bigr
\rrvert. %
\]
%
It is always dominated by the Kantorovich--Rubinstein distance since
\begin{eqnarray*}
\dTVRandom(\zeta_1,\zeta_2) & = &\sup
_{A\in\cN_\YY}\biggl\llvert\inf_{\bC\in\Sigma(\zeta_1,\zeta_2)} \int
_{\N_\YY\times\N_\YY} {\mathbf1}(\omega_1\in A)-{\mathbf1}(
\omega_2\in A) \bC\bigl(\dint(\omega_1,
\omega_2)\bigr)\biggr\rrvert
\\
& \leq&\inf_{\bC\in\Sigma(\zeta_1,\zeta_2)} \int
_{\N_\YY\times
\N_\YY} \dTV(\omega_1,\omega_2) \bC
\bigl(\dint(\omega_1,\omega_2)\bigr) = \dKR(
\zeta_1,\zeta_2) .
\end{eqnarray*}
The following example shows that convergence in Kantorovich--Rubinstein
distance is strictly finer than convergence in total variation distance.

%
\begin{example}
Let $\zeta$ be a Poisson process on $\YY$ with finite intensity
measure $\bM$. Let $(X_i)_{i\in\NN}$ be a sequence of independent
random elements in $\YY$ with distribution $\bM(\YY)^{-1}\bM(\cdot
)$ and let $Z$ be a Bernoulli random variable such that $\bP(Z=1)=p$
for some $p\in(0,1)$. Moreover, assume that $\zeta$, $(X_i)_{i\in\NN
}$ and $Z$ are independent. Now we consider the point process
\[
\zeta_{n,p}:=\zeta+{\mathbf1}(Z=1) \sum_{i=1}^n
\delta_{X_i} . %
\]
Since $\zeta$ and $\zeta_{n,p}$ coincide on an event with probability
$1-p$, we have that $\dTVRandom(\zeta,\zeta_{n,p}) \leq p$. By
taking $h(\mu)=\mu(\YY)$ as a test function in \eqref
{eqRubinsteinDistanceEquivalent}, we deduce that $\dKR(\zeta,\zeta
_{n,p}) \geq np$. Taking $p_n=1/\sqrt{n}$ for $p$ shows that
\[
\dTVRandom(\zeta,\zeta_{n,p_n})\to0 \quad\mbox{and}\quad\dKR(\zeta,
\zeta_{n,p_n}) \to\infty\qquad\mbox{as } n\to\infty, %
\]
so that $(\zeta_{n,p_n})_{n\in\NN}$ converges to $\zeta$ in total
variation distance but not in Kantorovich--Rubinstein distance.
\end{example}

In the previous example the Kantorovich--Rubinstein distance is more
strongly affected by the rare event that $\zeta\neq\zeta_{n,p_n}$
than the total variation distance, since the class of test functions is
larger and contains functions taking also values different from zero
and one. As already mentioned in the \hyperref[sec1]{Introduction}, one can say that the
difference between the Kantorovich--Rubinstein distance and the total
variation distance for point processes is similar to the difference
between the Wasserstein and the total variation distance for
integer-valued random variables. As particular example we cite the work
of Barbour and Xia \cite{BarbourXia2006}, where Poisson approximation
of random variables with respect to the Wasserstein distance has been
considered, extending previous results for the total variation
distance; see also Section~\ref{secUstatistics} below.

\subsection{A discrete gradient}\label{sec2.6}
For a counting measure $\omega\in\widetilde{\N}_\YY$
and a measurable
function $h\dvtx \widetilde{\N}_\YY\to{\mathbb{R}}$, let us introduce
the discrete
gradient in direction $y\in\YY$ by
\[
D_yh(\o):=h(\o+\d_y)-h(\o) , %
\]
where we recall that $\delta_y$ is the unit-mass Dirac measure
charging $y\in\YY$. In our notation we often suppress the dependence
of $D_yh(\o)$ on the underlying counting measure $\o$ and write
$D_yh$. Clearly, if $h\in\cC$, it holds that $\llvert D_yh\rrvert\leq
1$ for all
\mbox{$y\in\YY$}.

\subsection{Geometric preparations} For our applications in Section~\ref{secApplications},
we need some facts from convex
geometry. The Euclidean
norm in $\RR^d$ is denoted by $\llVert\cdot\rrVert$. The Euclidean distance
between two sets $A_1,A_2\subset\RR^d$ is given by
\[
\dist(A_1,A_2)=\inf\bigl\{\llVert x_1-x_2
\rrVert\dvtx x_1\in A_1, x_2\in A_2
\bigr\} . %
\]
If $A_1=\{x\}$ with $x\in\RR^d$, we write $\dist(x,A_2)$ instead of
$\dist(\{x\},A_2)$.
For a measurable set $K\subset\RR^d$, we write $\vol(K)$ for the volume
(i.e., $d$-dimensional Lebesgue measure) of $K$. For the volume of the
unit ball $B^d=\{x\in\RR^d\dvtx\llVert x\rrVert\leq1\}$ in $\RR^d$, we
introduce the
abbreviation $\kappa_d:=\vol(B^d)$. More generally, $B^d(x,r)$ will
denote the closed $d$-dimensional ball of radius $r>0$ centered at
$x\in\RR^d$, and we write $B^d(r)$ instead of $B^d(0,r)$ for short.
For $r\geq0$, the Minkowski sum $K_r=K+rB^d$ of $K$ and
$rB^d$ is the so-called $r$-parallel set of $K$. In particular, if $K$
is a convex set with nonempty interior, Steiner's formula (see, e.g.,
\cite{SW}, equation (14.5)) says
that the volume $\vol(K_r)$ is a polynomial of degree $d$ in
$r$. Formally,
%
\begin{equation}
\label{eqSteiner} \vol(K_r)=\sum_{i=0}^d
\kappa_{d-i} V_i(K) r^{d-i} .
\end{equation}
The coefficients $V_0(K),\ldots,V_d(K)$ are the so-called intrinsic
volumes of $K$, especially $V_0(K)=1$ whenever $K\neq\varnothing$,
$V_1(K)$ is a constant multiple of the mean width of $K$, $V_{d-1}(K)$
is half of the surface area of $K$ (if $K$ is the closure of its
interior) and $V_d(K)=\vol(K)$; cf. \cite{SW}, Chapter~14.2.

For $1\leq m\leq d-1$, we denote in this paper by $\GG_m^d$ the space
of $m$-dimensional linear subspaces and by $\AA_m^d$ the
space of $m$-dimensional affine subspaces of $\RR^d$. For
$L,M\in\GG_m^d$ let $[L,M]$ be the subspace determinant of $L$ and
$M$, that is, the $2m$-volume of a parallelepiped spanned by two
orthonormal bases in $L$ and in $M$. In one of our examples, we will also
deal with the integrated subspace determinant, and for this reason
we recall that
%
\begin{equation}
\label{eqIntegralL,M} \int_{\GG_m^d}\int_{\GG_m^d}[L,M]
\,\dint L \,\dint M=\frac
{{d-m\choose m}}{{d\choose m}}{\frac{\kappa_{d-m}^2}{\kappa_d\kappa_{d-2m}}}
\end{equation}
from \cite{HSS}, Lemma 4.4. Here, $\dint L$ and $\dint M$ indicate
integration with respect to the unique Haar probability measure on
$\GG_m^d$.

\section{Main results}\label{secmain-result}

\subsection{General estimate}\label{subsecGeneralEstimate}

Let $(\YY,\cY)$ be a lcscH space, and let
us fix another lcscH space $(\XX,\cX)$. We adopt the notation
introduced in Section~\ref{secPreliminaries} and denote by $\N_\XX$
the space of
$\sigma$-finite counting measures on $\XX$.

Let $\mu$ be a point process on $\XX$ with a $\sigma$-finite
intensity measure $\bK( \cdot):=\bE\mu( \cdot)$. Fix an
integer $k\geq1$, and let $f\dvtx\dom f \to\YY$ be a symmetric and measurable
function, where $\dom f$ is a symmetric subset of $\XX^k$; that is, if
$(x_1,\ldots,x_k)\in\dom f$, then
$(x_{\sigma(1)},\ldots,x_{\sigma(k)})\in\dom f$ for all permutations
$\sigma$ of $\{1,\ldots,k\}$. We now apply $f$ to all $k$-tuples of
distinct points of
$\mu$ contained in $\dom f$ to form a point process~$\xi$, that is,
\[
\xi:=\frac{1}{k!}\sum_{(x_1,\ldots,x_k)\in\mu_{\neq}^k\cap\dom
f}
\delta_{f(x_1,\ldots,x_k)} . %
\]
Since $f$ is symmetric, every $f(x_1,\ldots,x_k)$ also appears for the
$k!$ permutations of the argument $(x_1,\ldots,x_k)$. However, for
each subset $\{x_1,\ldots,x_k\}\subset\mu$ of distinct points of
$\mu$, we assign to
$f(x_1,\ldots,x_k)$ only multiplicity one as can be seen from the
above definition of $\xi$. However, $\xi$ might still have points of
multiplicity greater than one if there are different combinations of
$k$ points in $\XX$ that are mapped under $f$ to the same point
in $\YY$. The intensity measure of $\xi$ is denoted by $\bL$ and is
given by
\begin{eqnarray*}
\bL(A) &=& \bE\xi(A) = \bE\sum_{y\in\xi}{\mathbf1}(y\in A)
\\
&=&
\frac
{1}{k!} \bE\sum_{(x_1,\ldots,x_k)\in\mu^k_{\neq}\cap\dom f}{\mathbf1}
\bigl(f(x_1,\ldots,x_k)\in A\bigr) , \qquad A\in\cY.
\end{eqnarray*}

In what follows, we consider for $\mu$ two different types of point
processes, namely Poisson processes and binomial processes. By $\eta$
we denote a Poisson process on $\XX$ with a
$\sigma$-finite intensity measure $\bK$. By $\beta_n$ we
denote a binomial process of $n\in\NN$ points in $\XX$, which
are independent and identically distributed in $\XX$ according to a
probability measure $\bK_1$ on $\XX$. Such a binomial process
$\beta_n$ has intensity measure $\bK:=n\bK_1$. Now the multivariate
Mecke formula \eqref{eqMecke} and its binomial analogue \eqref
{eqMomentMeasureBinomial} imply that the intensity measure $\bL$ of
$\xi$ is given by
%
\begin{equation}
\label{eqLPoisson} \bL(A) =\frac{1}{k!}\int_{\dom f} {\mathbf1}
\bigl(f(x_1,\ldots,x_k)\in A\bigr) \bK^k\bigl(
\dint(x_1,\ldots,x_k)\bigr) , \qquad A\in\cY,\hspace*{-20pt}
\end{equation}
in the Poisson case and by
%
\begin{equation}
\label{eqLbinomial} \bL(A) =\frac{(n)_k}{k!} \int_{\dom f} {\mathbf1}
\bigl(f(x_1,\ldots,x_k)\in A\bigr) \bK_1^k
\bigl(\dint(x_1,\ldots,x_k)\bigr) , \qquad A\in\cY,\hspace*{-30pt}
\end{equation}
if we start with a binomial process. (To deal with both cases
simultaneously we use the same notation for both set-ups.) Let us
finally introduce $r(\dom f)$ for $k\geq2$ by
\begin{eqnarray*}
&& r(\dom f)
\\
&&\qquad  :=\max_{1\leq\ell\leq k- 1} \int_{\XX^\ell} \biggl(
\int_{\XX^{k-\ell}}{\mathbf1}\bigl((x_1,\ldots,x_k)
\in\dom f\bigr) \bK^{k-\ell}\bigl(\dint(x_{\ell+1},
\ldots,x_{k})\bigr) \biggr)^2
\\
&&\hspace*{86pt}{}\times \bK^\ell\bigl(
\dint(x_1,\ldots,x_\ell)\bigr) , %
\end{eqnarray*}
and, for $k=1$, put $r(\dom f):=0$. Moreover, we use the convention
that $(n-k)_k/(n)_k:=0$ if $n<k$.

We can now state our main result, a functional limit theorem, which
provides a bound on the Kantorovich--Rubinstein distance between $\xi$
and a suitable Poisson process on $\YY$.

%
\begin{theorem}\label{thmMain}
Let $\zeta$ be a Poisson process on $\YY$ with finite intensity
measure~$\bM$. If $\xi$ is induced by the Poisson process $\eta$, then
%
\begin{eqnarray}
\label{eqPPPconv6} \dKR(\xi,\zeta) & \leq&\dTV(\bL,\bM)+2 \bigl(\bE\xi(
\YY)^2 - \bE\xi(\YY)-\bigl(\bE\xi(\YY)\bigr)^2 \bigr)
\nonumber\\[-8pt]\\[-8pt]\nonumber
& \leq&\dTV(\bL,\bM)+ \frac{2^{k+1}}{k!} r(\dom f) .
\end{eqnarray}
If otherwise $\xi$ is derived from the binomial process $\beta_n$, then
\begin{eqnarray*}
\dKR(\xi,\zeta) & \leq&\dTV(\bL,\bM)+2 \biggl(\bE\xi(
\YY)^2-\bE\xi(\YY)-\frac{(n-k)_k}{(n)_k}\bigl(\bE\xi(\YY)
\bigr)^2 \biggr)
\\
&&{} +\frac{6^k
k!}{n} \bigl(\bE\xi(\YY)
\bigr)^2
\\
& \leq&\dTV(\bL,\bM)+\frac{2^{k+1}}{k!}r(\dom f)+\frac{6^k k!}{n} \bL(
\YY)^2 .
\end{eqnarray*}
\end{theorem}

%
\begin{remark}\label{remMain}
(i)~If the underlying point process is a binomial process
$\beta_n$ with $n$ points and if $n<k$, the point process $\xi$ is
empty with probability one and $\bL\equiv0$. In this case, $\dKR
(\xi,\zeta)\leq\bE\zeta(\YY)=\dTV(\bL,\bM)$, and the bound on
$\dKR(\xi,\zeta)$ is trivially valid. For this reason, no
further restriction on $n$ is necessary.

(ii)~In the Poisson case, it can happen that $\bL(\YY)=\infty
$. In this case, we have $\dTV(\bL,\bM)=\infty$, and the bound
\eqref{eqPPPconv6} is trivial. Hence Theorem~\ref{thmMain} is only
of interest if $\bL(\YY)<\infty$, which is equivalent to $\bK
^k(\dom f)<\infty$, a condition which ensures that $\xi$ is almost
surely finite.

(iii)~Taking $\bM=\bL$ in the Poisson case in Theorem~\ref
{thmMain} shows that
\[
\dKR(\xi,\zeta)\leq2 \bigl(\bE\xi(\YY)^2 - \bE\xi(\YY)-\bigl(\bE
\xi(\YY)\bigr)^2 \bigr)\leq\frac{2^{k+1}}{k!}r(\dom f) . %
\]
In particular, if $k=1$, this gives $\dKR(\xi,\zeta)=0$, which in
view of Proposition~\ref{thmPPPconv1} implies that $\xi$ is a
Poisson process. This is consistent with the well-known mapping theorem
for Poisson processes, for which we refer to \cite{Kingman}, Chapter~2.3.

(iv) If $\XX=\YY$ and $f\dvtx\XX\to\XX$ is the identity,
Theorem~\ref{thmMain} yields that, for Poisson processes $\xi$ and
$\zeta$ with finite intensity measures $\bL$ and $\bM$, respectively,
\[
\dKR(\xi,\zeta)\leq\dTV(\bL,\bM) . %
\]
In other words, the Kantorovich--Rubinstein distance between two
Poisson processes is bounded by the total variation distance of their
intensity measures. For a similar estimate in a more restricted set-up
we refer to \cite{DecreusefondRubinstein}, Proposition 4.1.
\end{remark}

\subsection{The Euclidean case}\label{secSpecialCases}

In this subsection we shall apply our general estimate of Theorem
\ref{thmMain} to the important situation that the target space $\YY$
is $\RR^d$ endowed with the standard
Borel $\sigma$-field $\mathcal{B}(\RR^d)$. This is tailored toward
some of our
applications in Section~\ref{secApplications} and is similar to the
set-up in \cite{STScalingLimits,STFlats}. We let
$(\XX,\cX)$ be a lcscH space and let $(\eta_t)_{t\geq1}$ be a
family of Poisson processes in $\XX$ with intensity measures $\bK
_t=t\bK$, $t\geq1$, where $\bK$ is a fixed $\sigma$-finite measure.
By $(\beta_t)_{t\geq1}$ we denote a family of binomial processes such
that $\beta_t=\beta_{\lceil t\rceil}$, and $\beta_{\lceil t\rceil
}$ is a process of $\lceil t \rceil$ points chosen independently
according to a fixed probability measure $\bK_1$. In this situation we
use the notation $\bK_t:=\lceil t\rceil\bK$. We write $(\mu
_t)_{t\geq1}$ in the sequel to indicate either $(\eta_t)_{t\geq1}$
or $(\beta_t)_{t\geq1}$.

For a fixed integer $k\geq1$ we consider symmetric and measurable
functions $f_t\dvtx \XX^k \to\RR^d$, $t\geq1$. We are interested in the
behavior of the derived point processes
\[
\xi_t:=\frac{1}{k!} \sum_{(x_1,\ldots,x_k)\in\mu_{t,\neq
}^k}
\delta_{f_t(x_1,\ldots,x_k)} , \qquad t\geq1 . %
\]
For this reason, we consider the re-scaled point processes
\[
t^\g\bullet\xi_t:=\frac{1}{k!} \sum
_{(x_1,\ldots,x_k)\in\mu
_{t,\neq}^k}\d_{t^\g
f_t(x_1,\ldots,x_k)} , \qquad t\geq1 , %
\]
where $\gamma\in\RR$ is a suitable constant. In order to compare
$t^\g\bullet\xi_t$ with a Poisson process, we need to introduce the following
notation. The intensity measure $\bL_t$ of the re-scaled point process
$t^\g\bullet\xi_t$ is given by
\[
\bL_t(B):=\frac{1}{k!} \bE\sum_{(x_1,\ldots,x_k)\in\mu_{t,\neq
}^k}{
\mathbf1} \bigl(f_t(x_1,\ldots,x_k)\in
t^{-\g}B \bigr) ,\qquad B\in{\mathcal B}\bigl(\RR^d\bigr) .
\]
For $B\in\mathcal{B}(\RR^d)$ let $r_t(B)$ be given by $r_t(B):=0$
for $k=1$ and
\begin{eqnarray*}
r_t(B) &:=& \max_{1\leq\ell\leq k- 1} \int_{\XX^\ell}
\biggl(\int_{\XX^{k-\ell}}{\mathbf1}\bigl(f_t(x_1,
\ldots,x_k) \in t^{-\g}B\bigr) \bK_t^{k-\ell}
\bigl(\dint(x_{\ell+1},\ldots,x_k)\bigr) \biggr)^2
\\
&&\hspace*{51pt}{}\times \bK_t^\ell\bigl(\dint(x_1,\ldots,x_\ell)
\bigr) %
\end{eqnarray*}
for $k\geq2$. Furthermore, for a measure $\nu$ on $\RR^d$ and $B\in
{\mathcal B}(\RR^d)$ let $\nu\mid_B$ be the restriction of $\nu$ to $B$.

%
\begin{corollary}\label{corScalingLimits}
Let $\zeta$ be a Poisson process on $\RR^d$ with intensity measure~$\bM$, and let $B\in\mathcal{B}(\RR^d)$ be such that $\bM
(B)<\infty$. If $\xi_t$ is induced by a Poisson process $\eta_t$
with $t\geq1$, then
\begin{eqnarray*}
&& \dKR\bigl(\bigl(t^\g\bullet\xi_t\bigr)\llvert
_B,\zeta\rrvert_B \bigr)
\\
&&\qquad  \leq \dTV\bigl(
\bL_t\llvert_B,\bM\rrvert_B\bigr)+2
\bigl(\bE\xi_t\bigl(t^{-\gamma}B\bigr)^2-\bE\xi
_t\bigl(t^{-\gamma}B\bigr)-\bigl(\bE\xi_t
\bigl(t^{-\gamma}B\bigr)\bigr)^2 \bigr)
\\
&&\qquad  \leq \dTV\bigl(\bL_t\llvert_B,\bM\rrvert
_B\bigr)+\frac{2^{k+1}}{k!}r_t(B) .
\end{eqnarray*}
If $\xi_t$ is induced by a binomial process $\b_t$ with $t\geq1$, then
\begin{eqnarray*}
&& \dKR\bigl(\bigl(t^\g\bullet\xi_t\bigr)\llvert
_B,\zeta\rrvert_B \bigr)
\\
&&\qquad  \leq\dTV\bigl(
\bL_t\llvert_B,\bM\rrvert_B\bigr)
\\
&&\quad\qquad{} +2
\biggl(\bE\xi_t\bigl(t^{-\gamma}B\bigr)^2-\bE\xi
_t\bigl(t^{-\gamma}B\bigr)-\frac{(\lceil t\rceil-k)_k}{(\lceil t\rceil
)_k}\bigl(\bE
\xi_t\bigl(t^{-\gamma}B\bigr)\bigr)^2 \biggr)
\\
&&\quad\qquad{} +\frac{6^k
k!}{t} \bigl(\bE
\xi_t\bigl(t^{-\gamma}B\bigr)\bigr)^2
\\
&&\qquad \leq \dTV\bigl(\bL_t\llvert_B,\bM\rrvert
_B\bigr)+\frac{2^{k+1}}{k!}r_t(B)+\frac{6^k
k!}{t}
\bL_t(B)^2 .
\end{eqnarray*}
\end{corollary}
\begin{pf}
This is a direct consequence of Theorem~\ref{thmMain} with $t^\gamma
\bullet\xi_t\mid_B$ instead of $\xi$ and $\zeta\mid_B$ instead of
$\zeta
$ there.
\end{pf}

In view of limit theorems, the most natural choice for $\bM$ is to
take $\bM$ as the strong limit of the measures $\bL_t$, as $t\to
\infty$. That is,
\[
\bM(B)=\lim_{t\to\infty}\bL_t(B)\qquad\mbox{for all }B
\in\cB\bigl(\RR^d\bigr) .
\]
However, we emphasize that this does not necessarily imply that $\dTV
(\bL_t,\break \bM)\to0$, as $t\to\infty$, even though this is true for
our applications presented below.

%
\begin{remark}\label{remnonUniform}
(i) The upper bounds in Corollary~\ref{corScalingLimits} are not
uniform in the sense that they depend on the set $B$. This was to be
expected since the re-scaled point processes $t^{\g}\bullet\xi_t$
can be
finite for any $t\geq1$, while a realization of $\zeta$ can charge an
infinite number of
points (compare with our applications in Section~\ref
{secApplications}). This is the reason for introducing the restriction
to the
set $B$, which allows us to compare $t^{\g}\bullet\xi_t\mid_B$ with
$\zeta\mid_B$ using the Kantorovich--Rubinstein distance.

(ii) To allow for an easier comparison with the previous
paper \cite{STScalingLimits}, we remark that ibidem the Poisson case
for $d=1$ is considered. Moreover, the intensity measure $\bM$
there is concentrated on the positive real half-axis and has the form
\[
\bM(B)=a b \int_B {\mathbf1}(u\geq0) u^{b-1} \,\dint
u , \qquad B\in\mathcal{B}(\RR) , %
\]
for some constants $a,b> 0$. In this case, the
Poisson process $\zeta$ is a so-called Weibull process since the
distance from the origin to the closest point of $\zeta$ is Weibull
distributed with distribution function $u\mapsto(1-\exp(-a u^b))
{\mathbf1}(u>0)$. We remark that this form of $\bM$
was tailored to the applications in \cite{STScalingLimits}; a more
general version is stated without proof in \cite{STFlats}.

(iii) Note that $r_t(B)$ is dominated by $k!\bL_t(B)\hat
r_t(B)$, where $\hat{r}_t(B)$ is defined as
\begin{eqnarray*}
\hat r_t(B)
&:=& \mathop{\max_{1\leq\ell\leq k- 1,}}_{(x_1,\ldots
,x_\ell)\in\XX^\ell}
\bK_t^{k-\ell} \bigl(\bigl\{(y_1,
\ldots,y_{k-\ell})\in\XX^{k-\ell
}\dvtx
\\[-1.5pt]
&&\hspace*{84pt}  f_t(x_1,
\ldots,x_{\ell},y_1,\ldots,y_{k-\ell})\in
t^{-\g}B\bigr\} \bigr) %
\end{eqnarray*}
for $B\in\cB(\RR^d)$. A quantity similar to $\hat r_t(B)$ has also
played a prominent role in the previous study
\cite{STScalingLimits}. In many applications a bound for $\hat
r_t(B)$ is already sufficient in order to apply Corollary
\ref{corScalingLimits}. However, there are situations for which
$\hat r_t(B)$ is an increasing function in $t$, while $r_t(B)$
tends to zero, as $t\to\infty$. This way \cite{STScalingLimits},
Theorem~1.1, in which $\hat r_t$ instead of $r_t$ appears,
is not applicable in such cases, as is erroneously done in Sections~\ref{sec2.5} and
\ref{sec2.6} ibidem. However, in these specific cases it is readily
checked that $r_t$ behaves nicely, implying that the results there are correct.
\end{remark}

\section{A general Stein principle}\label{secgener-stein-princ}

This section is devoted to a more informal discussion about the method
of bounding the Kantorovich--Rubinstein distance between point
processes using a Stein principle. This approach is the key argument of
our proof of Theorem~\ref{thmMain} in Section~\ref
{secProof-main-result}. Recall that the aim is to provide an upper
bound for the Kantorovich--Rubinstein distance between a Poisson
process $\zeta$ on a space $\YY$ with finite intensity measure $\bM$
and a second point process $\xi$ on $\YY$, which in turn is derived
from another point process $\mu$ on a space $\XX$ by a transformation.

The first part of Stein's method consists of characterizing the target
object, here the Poisson process $\zeta$. The method is to consider a
functional operator $L$ which, at a formal level, satisfies for a
finite point process $\nu$
the identity
%
\begin{equation}
\label{eqPPPconv1} \bE\bigl[LF(\nu)\bigr] =0\qquad\mbox{for a large
class of
functions $F\dvtx\widehat{\N}_\YY\to\RR$}
\end{equation}
if and only if $\nu$ is a Poisson process with intensity measure
$\bM$. It is usually not difficult to construct such an operator for a
given target object. What may become challenging, especially in
infinite dimensions (compare with
\cite{Barbour90,DecreusefondCoutinInfiniteStein,Shih2011}), will be to prove
that the target object is the unique solution of \eqref{eqPPPconv1}.
In our case, uniqueness follows from the theory of spatial birth--death
processes; see \cite{Preston75}.

The second step of Stein's method is to solve the so-called Stein
equation,
%
\begin{equation}
\label{eqStein1} LF(\omega)=\bE h(\zeta) - h(\omega) , \qquad\omega\in
\widehat{
\N}_\YY,
\end{equation}
for a certain class of test functions $h\dvtx\widehat{\N}_\YY\to\RR$.
This means that we have to compute a solution $F_h$ for a given test
function $h$ and to evaluate $LF_h(\omega)$.

A prominent way to do this is to use the so-called generator approach;
see the survey article \cite{Reinert2005} and the references cited
therein. The underlying idea is to interpret $L$ as infinitesimal
generator of a Markov process with the distribution of $\zeta$ as its
invariant distribution, whence $L$ satisfies \eqref{eqPPPconv1}. If
$(P_s)_{s\geq0}$ is the semi-group associated with this Markov
process, one can show that
%
\begin{equation}
\label{eqStein2} L F_h(\omega)= \int_0^\infty
LP_sh(\omega) \,\dint s , \qquad\omega\in\widehat{
\N}_\YY.
\end{equation}
In order to compare the point process $\xi$ with $\zeta$, we put
$\omega=\xi$ and take expectations in \eqref{eqStein1} and \eqref
{eqStein2}. This leads to
\[
\bE h(\zeta) - \bE h(\xi) = \bE L F_h = \bE\int
_0^\infty LP_sh(\xi) \,\dint s .
\]
In the subsequent section, we will derive this identity rigorously. In
the context of our main result, the point process $\xi$ is induced by
an underlying point process $\mu$ on another space $\XX$. More
formally we have that $\xi=T(\mu)$, where $T$ is a suitable
transformation, that is, a mapping from $\N_\XX$ to $\widehat{\N
}_\YY$. Hence we will have to compute
\[
\bE\int_0^\infty LP_sh\bigl(T(\mu)
\bigr) \,\dint s . %
\]
This expression is bounded in Section~\ref{secProof-main-result} by
exploiting the special structure of the transformation $T$ and the fact
that $\mu$ is a Poisson or binomial process.

\section{Glauber dynamics for the Poisson process} \label{secglaub-dynam-poiss}

We now specialize the general scheme outlined in Section~\ref
{secgener-stein-princ} to our particular situation. Although the
approach is similar to \cite{BarbourBrown}, Section~2, for example, we
prefer to carry out the details here since we consider a different
class of test functions, namely Lipschitz functions instead of bounded
functions. We assume the same set-up as for Theorem~\ref{thmMain};
that is, $\zeta$ is a Poisson process on a lcscH space $\YY$ with a
finite intensity measure $\bM$ and distribution $\bP_\zeta$. We now
construct a Glauber dynamics
for $\bP_\zeta$, that is a continuous-time Markov process
$(G(s))_{s\geq0}$ with state space $\widehat{\N}_\YY$ and $\bP
_\zeta$ as its stationary (i.e., invariant) distribution; see
\cite{Preston75}. Its generator $L$ is given by
%
\begin{eqnarray}\label{eqGeneratorL}
\quad Lh(\o):=\int_\YY h(\o+\d_y)-h(
\o) \bM(\dint y)+\int_\YY h(\o-\d_y)-h(\o) \o(
\dint y) ,
\nonumber\\[-8pt]\\[-8pt]
\eqntext{\o\in\widehat{\N}_\YY,}
\end{eqnarray}
where $h\dvtx\widehat{\N}_\YY\to\RR$ is a measurable and bounded
function. According to our notational convention, $L$ may be re-written
as
\[
Lh(\o)=\int_\YY h(\o+\d_y)-h(\o) \bM(\dint
y)+\sum_{y\in\o} \bigl(h(\o-\d_y)-h(\o)
\bigr) . %
\]
Note that $Lh(\o)$ is well defined for all $h\in\cC$ and $\o\in
\widehat{\N}_\YY$ since the Lipschitz property implies that the
integrands in \eqref{eqGeneratorL} are bounded by one. Moreover, we
notice that the operator $L$ uniquely determines the process
$(G(s))_{s\geq0}$, which has $\bP_\zeta$ as its unique invariant
distribution; see \cite{DaleyVereJones2008}, Proposition 10.4.VII, or
\cite{Preston75}.

The Markov process $(G(s))_{s\geq0}$ is a spatial birth--death process in
continuous time whose dynamics can be described as follows. If at time
$s$, the system is in state~$\o_s$, each particle charged by $\o_s$
dies at rate $1$, and a new particle is born at $y$ with rate
$\bM(\dint y)$. Alternatively, imagine a homogeneous Poisson process
$\zeta_b$ on $\RR_+$ with intensity $\bM(\YY)$. The jump times of
$\zeta_b$ determine the birth times of the particles in $\zeta$. At
each jump of $\zeta_b$ a new particle is born and is placed in $\YY$
according to the distribution $\bM( \cdot)/\bM(\YY)$,
independently of the current configuration. Moreover, each particle has
a lifetime which is exponentially distributed with parameter~$1$,
independent of the past and of the rest of the configuration; see
again~\cite{Preston75}.

The semi-group $(P_s)_{s\geq0}$ associated with the Markov process
$(G(s))_{s\geq0}$ is defined as
%
\begin{equation}
\label{eqSemigroupPt} P_sh(\o)=\bE\bigl[ h\bigl(G(s)\bigr) \mid G(0)=\o
\bigr] , \qquad\o\in\widehat{\N}_\YY,  h\dvtx\widetilde{
\N}_\YY\to\RR.
\end{equation}
For $h\in\cC$ and $\omega\in\widehat{\N}_\YY$ the conditional
expectation is always well defined since
\begin{eqnarray*}
\bigl\llvert P_sh(\o)\bigr\rrvert&=&\bigl\llvert\bE\bigl[ h
\bigl(G(s)\bigr) \mid G(0)=\o\bigr]\bigr\rrvert
\\
& \leq&\bE\bigl[ \bigl\llvert h\bigl(G(s)\bigr)-h(\o) \bigr\rrvert\mid
G(0)=\o
\bigr]+\bigl\llvert h(\o)\bigr\rrvert
\\
& \leq&\bE\bigl[ \dTV\bigl(G(s),\o\bigr) \mid G(0)=\o\bigr]+\bigl
\llvert h(\o)
\bigr\rrvert
\\
& \leq&\bE\zeta_b\bigl([0,s]\bigr) + \o(\YY) + \bigl\llvert h(\o)
\bigr\rrvert<\infty,
\end{eqnarray*}
where $\zeta_b$ is the homogeneous Poisson process from the
description of the birth--death dynamics above. Below we will need the
following lemmas
about the process $(G(s))_{s\geq0}$ and its semi-group $(P_s)_{s\geq
0}$. The first one provides a commutation relation between the discrete
gradient and the semi-group.

%
\begin{lemma}\label{lemDPt}
For any $s\geq0$, $\omega\in\widehat{\N}_\YY$, $y\in\YY$ and
$h\in\cC$,
\[
D_yP_sh(\o)=e^{-s} P_s(D_yh)
(\o) . %
\]
\end{lemma}
\begin{pf}
To construct a sample path of $(G(s))_{s\geq0}$, given the initial
configuration $G(0)=\o+\d_y$, we have to add the
independent particle $y$ to a realization of $(G(s))_{s\geq0}$
starting from the initial configuration $\o$. These two realizations
will be identical after the particle $y$ has died. Thus, denoting
by $\ell(y)$ the lifetime of $y$ and using
\eqref{eqSemigroupPt}, we can write
\begin{eqnarray*}
D_yP_sh(\o) &=& \bE\bigl[h\bigl(G(s)\bigr) \mid G(0)=
\o+\d_y\bigr]-\bE\bigl[h\bigl(G(s)\bigr) \mid G(0)=\o\bigr]
\\
&=& \bE\bigl[\bigl(h\bigl(G(s)+\d_y\bigr)-h\bigl(G(s)\bigr)\bigr) {
\mathbf1}\bigl(\ell(y)\geq s\bigr) \mid G(0)=\o\bigr] .
\end{eqnarray*}
Since $\ell(y)$ is independent of everything else and is
exponentially distributed with mean one, we can continue with
\begin{eqnarray*}
D_yP_sh(\o) &=& \bE\bigl[{\mathbf1}\bigl(\ell(y)\geq s
\bigr)\bigr] \bE\bigl[\bigl(h\bigl(G(s)+\d_y\bigr)-h\bigl(G(s)\bigr)
\bigr) \mid G(0)=\o\bigr]
\\
&=& e^{-s} P_s(D_yh)(\o),
\end{eqnarray*}
where we have used \eqref{eqSemigroupPt} again. This completes the
proof.
\end{pf}

%
\begin{lemma}\label{lemPreparationForOtherLemmas}
Let $\o_1,\o_2\in\widehat{\N}_\YY$ with $\o_2\subset\o_1$. If
$h\in\cC$ and $s\geq0$, then
\[
\bigl\llvert\bE\bigl[h\bigl(G(s)\bigr) \mid G(0)=\o_1\bigr]-\bE
\bigl[h\bigl(G(s)\bigr) \mid G(0)=\o_2\bigr]\bigr\rrvert\leq(
\o_1\setminus\o_2) (\YY) e^{-s} . %
\]
\end{lemma}
\begin{pf}
Recall that each particle $y$ of the initial configuration $G(0)$
has an exponentially distributed lifetime $\ell(y)$ with mean
one. Thus since $h\in\cC$, it holds that
\begin{eqnarray*}
&& \bigl\llvert\bE\bigl[h\bigl(G(s)\bigr) \mid G(0)=\o_1\bigr]-\bE
\bigl[h\bigl(G(s)\bigr) \mid G(0)=\o_2\bigr]\bigr\rrvert
\\
&&\qquad\leq\bE\biggl[\biggl\llvert h\biggl(G(s)+\sum
_{y\in\o_1\setminus\o
_2}{\mathbf1}\bigl(\ell(y)\geq s\bigr)\d_y
\biggr)-h\bigl(G(s)\bigr)\biggr\rrvert\Big| G(0)=\o_2 \biggr]
\\
&&\qquad\leq\bE\biggl[\dTV\biggl(G(s)+\sum_{y\in\o_1\setminus\o_2}{
\mathbf1}\bigl(\ell(y)\geq s\bigr)\d_y,G(s)\biggr) \Big| G(0)=
\o_2 \biggr]
\\
&&\qquad\leq\bE\sum_{y\in\o_1\setminus\o_2}{\mathbf1}\bigl(\ell(y)\geq s
\bigr)
\\
&&\qquad= (\o_1\setminus\o_2) (\YY) e^{-s} ,
\end{eqnarray*}
which proves the claim.
\end{pf}

%
\begin{lemma}\label{lemErgodic}
For any $\o\in\widehat{\N}_\YY$ and $h\in\cC$,
\[
\lim_{s\to\infty}P_sh(\o)=\bE h(\zeta)=\int h \,\dint
\bP_\zeta. %
\]
\end{lemma}
\begin{pf}
We notice first that the expectation on the right-hand side is well
defined since $h\in\cC$
implies that
\begin{eqnarray*}
\bE\bigl\llvert h(\zeta)\bigr\rrvert&\leq&\bE\bigl\llvert h(\zeta
)-h(\varnothing)
\bigr\rrvert+\bigl\llvert h(\varnothing)\bigr\rrvert\leq\bE\dTV(\zeta
,\varnothing)+
\bigl\llvert h(\varnothing)\bigr\rrvert\leq\bE\zeta(\YY)+\bigl\llvert
h(\varnothing)
\bigr\rrvert
\\
&=& \bM(\YY)+\bigl\llvert h(\varnothing)\bigr\rrvert, %
\end{eqnarray*}
where $\varnothing$ stands for the counting measure that corresponds
to the empty point configuration.

From Lemma~\ref{lemPreparationForOtherLemmas} with $\o_1=\o$ and
$\o_2=\varnothing$, we have that
%
\begin{equation}
\label{eqBoundCoupling1} \bigl\llvert\bE\bigl[h\bigl(G(s)\bigr) \mid
G(0)=\o\bigr] - \bE
\bigl[h\bigl(G(s)\bigr) \mid G(0)=\varnothing\bigr]\bigr\rrvert\leq
\omega(\YY)
e^{-s} .
\end{equation}
The number of particles of $G(s)$ starting from the
empty configuration follows the evolution of an M${}/{}$M${}/\infty$ queue
with arrival (birth) rate $\bM(\YY)$ and service (death) rate $1$,
and thus is Poisson
distributed with parameter $(1-e^{-s})\bM(\YY)$. Since the
position of each of the particles is independent of everything else,
$G(s)$ has the same distribution as a Poisson
process on $\YY$ with intensity measure $(1-e^{-s})\bM$. Since $\zeta
$ has the same distribution as the superposition of two independent
Poisson processes with intensity measures $(1-e^{-s})\bM$ and
$e^{-s}\bM$, respectively, we obtain that
%
\begin{equation}
\label{eqBoundCoupling2} \bigl\llvert\bE\bigl[h\bigl(G(s)\bigr) \mid
G(0)=\varnothing\bigr]-
\bE h(\zeta)\bigr\rrvert\leq e^{-s} \bM(\YY) .
\end{equation}
Combining \eqref{eqBoundCoupling1} and \eqref{eqBoundCoupling2} and
letting $s\to\infty$ completes the proof.
\end{pf}

The next lemma, which can be seen as an integration by parts formula,
is the key for the proof of Theorem~\ref{thmMain} given in
Section~\ref{secProof-main-result} below.

%
\begin{lemma} \label{lemPPPconv1}
If $h\in\cC$ and $\omega\in\widehat{\N}_\YY$, then
%
\begin{equation}
\label{eqSteinIdentity1} \bE h(\zeta)- h(\o)=\int_0^\infty
LP_sh(\o) \,\dint s .
\end{equation}
\end{lemma}
\begin{pf}
For an arbitrary $h\in\cC$ we define $h_n\dvtx\widehat{\N}_\YY\to\RR
$, $n\in\NN$ by
\[
h_n(\o)=\cases{ n, &\quad$h({\o})> n$,
\vspace*{3pt}\cr
h({\o}), &\quad$-n\leq
h({\o})\leq n$,
\vspace*{3pt}\cr
-n, &\quad$h({\o})<-n$.} %
\]
Clearly, each of the functions $h_n$ is bounded and belongs to $\cC$.
Thus the forward-backward equation stated as Theorem 12.22 in \cite
{KallenbergFoundations} implies that
%
\begin{equation}
\label{eqBoundLPshn} P_t h_n(\o)- h_n(\o) =\int
_0^t LP_sh_n(\o)\,
\dint s , \qquad t\geq0 .
\end{equation}
By construction, we have $h_n(\omega)\to h(\omega)$, as $n\to\infty
$. The dominated convergence theorem implies that $P_sh_n(\o)\to
P_sh(\o)$ and $LP_sh_n(\o)\to LP_sh(\o)$, as $n\to\infty$, for all
$s\geq0$. By \eqref{eqGeneratorL} and Lemma~\ref{lemDPt}, we have
that, for $g=h$ or $g=h_n$ and $s\geq0$,
%
\begin{eqnarray}
\label{eqBoundLPsh} \bigl\llvert LP_sg(\omega)\bigr\rrvert&\leq&\int
_\YY e^{-s} \bigl\llvert P_s(D_yg)
(\o)\bigr\rrvert\bM(\dint y)\nonumber
\\
&&{} +\int_\YY e^{-s}
\bigl\llvert P_s(D_yg) (\o-\delta_y)\bigr
\rrvert\o(\dint y)
\\
&\leq& e^{-s} \bigl(\bM(\YY)+\o(\YY)\bigr) .\nonumber
\end{eqnarray}
In the last step we used the fact that $\llvert P_s(D_yg)\rrvert\leq
1$. Now, a
further application of the dominated convergence theorem shows that
\[
\lim_{n\to\infty} \int_0^t
LP_sh_n(\o) \,\dint s = \int_0^t
LP_sh(\o) \,\dint s , \qquad t\geq0 , %
\]
so that, letting $n\to\infty$ in \eqref{eqBoundLPshn}, yields
%
\begin{equation}
\label{eqIntegralRepresentationPt} P_t h(\o)- h(\o) =\int_0^t
LP_sh(\o) \,\dint s , \qquad t\geq0 .
\end{equation}
Because of \eqref{eqBoundLPsh} and the dominated convergence theorem,
the right-hand side of \eqref{eqIntegralRepresentationPt} converges
to the right-hand side of \eqref{eqSteinIdentity1}, as $t\to\infty
$. Together with Lemma~\ref{lemErgodic} for the left-hand side, this
completes the proof.
\end{pf}

%
\begin{remark}
The operator $L$ and the associated semi-group $(P_s)_{s\geq0}$ on the
Poisson space can be also defined via the Wiener--It\^o chaos expansion,
which we recall now for completeness. We still denote by $\zeta$ a
Poisson process with intensity measure $\bM$ on a lcscH space $\YY$.
A crucial property
of $\zeta$ is that any square integrable functional $F\in
L^2(\bP_\zeta)$ of $\zeta$ can be written as
%
\begin{equation}
\label{eqChaosExpansion} F=\bE F+\sum_{n=1}^\infty
I_n(f_n)
\end{equation}
with
\[
f_n(y_1,\ldots,y_n)=\frac{1}{n!} \bE
D^n_{y_1,\ldots,y_n}F(\zeta) ,\qquad y_1,
\ldots,y_n\in\YY,  n\geq1 , %
\]
where $D^n:=D\circ D^{n-1}$ with $D^1:=D$ is the $n$th iteration of the
discrete gradient $D$ introduced in Section~\ref{secPreliminaries}, and
where $I_n(f_n)$ stands for the $n$-fold Wiener--It\^o integral of the
square integrable and symmetric function $f_n$ with respect to the signed
random measure $\zeta-\bM$. Moreover, the series in
\eqref{eqChaosExpansion} converges in $L^2(\bP_\zeta)$ and is called
the Wiener--It\^o chaos expansion of $F$; we refer to \cite
{LastPenrose} for further
details.
We can now define the Ornstein--Uhlenbeck generator $L$ on the Poisson
space by
\[
LF=-\sum_{n=1}^\infty nI_n(f_n)
, %
\]
whenever $F$ belongs to
$\dom L$; that is, $F$ is such that $\sum_{n=1}^\infty
n^2 n!\llVert f_n\rrVert_{L^2(\bM^n)}^2<\infty$, where
$\llVert\cdot\rrVert_{L^2(\bM^n)}$ stands for the usual norm in
$L^2(\bM^n)$. We\vspace*{1pt} remark that $LF$ can equivalently be written as in
\eqref{eqGeneratorL}
as a consequence of identity (3.19) in \cite{LastPenrose} and of the
relation stated in \cite{PSTU10}, Lemma 2.11, between the discrete
gradient, the Ornstein--Uhlenbeck generator and the so-called
Skorohod-integral on the Poisson space, another operator, which is not
needed in the sequel. In \cite{LastPeccatiSchulte2014} the relation
between the inverse of the Ornstein--Uhlenbeck generator and the
associated semi-group is investigated. The semi-group $(P_s)_{s\geq0}$
can be written in terms of the Wiener--It\^o chaos expansion as
\[
P_sF=\bE F +\sum_{n=1}^\infty
e^{-ns}I_n(f_n) ,\qquad s\geq0 , %
\]
where $F\in\dom L$ is assumed to have a chaotic expansion as in
\eqref{eqChaosExpansion}; see, for example, \cite
{LastPeccatiSchulte2014}, equation~(3.13). Lemma~\ref{lemDPt} is a
special case
of \cite{LastPeccatiSchulte2014}, Lemma 3.1, and Lemmas~\ref{lemPreparationForOtherLemmas}, \ref{lemErgodic} and~\ref{lemPPPconv1} can also be derived via the approach sketched in
this remark. However, we prefer to give proofs not relying on
Wiener--It\^o chaos expansions rather than on trajectorial properties.
\end{remark}

%
\begin{remark}
In \cite{SchuhmacherStucki} a spatial birth--death process is
constructed whose invariant distribution is a Gibbs process. This
includes the birth--death process in the present paper as a special
case, and the generator in \cite{SchuhmacherStucki} is a
generalization of the generator in \eqref{eqGeneratorL}. However, the
results in \cite{SchuhmacherStucki} do not cover the results of this
section since only the test functions for the total variation distance
are considered, while we use Lipschitz functions, which are needed for
the Kantorovich--Rubinstein distance.
\end{remark}

\section{Proof of Theorem \texorpdfstring{\protect\ref{thmMain}}{3.1}}\label{secProof-main-result}

Before going into the details of the proof of Theorem~\ref{thmMain},
we explain the strategy
informally in case of an underlying Poisson process $\eta$. Applying the
multivariate Mecke formula \eqref{eqMecke} in equation
\eqref{eqdifferenceExpectations} below,\vadjust{\goodbreak} we are lead to estimate the
integral with respect to $\bK^k$ of
%
\begin{eqnarray}\label{eqPPPconv5}
\qquad \bE\bigl[F\bigl(\xi(\eta+\d_{x_1}+\cdots+
\d_{x_k})-\d_{f(x_1,\ldots
,x_k)}\bigr) - F\bigl(\xi(\eta+\d_{x_1}+
\cdots+\d_{x_k})\bigr) \bigr] ,
\nonumber\\[-8pt]\\[-8pt]
\eqntext{x_1,\ldots,x_k\in\XX,}
\end{eqnarray}
with $F\dvtx\widehat{\N}_\YY\to\RR$ being a certain point process
functional and
where we write $\xi(\mu)$ instead of $\xi$ to underpin the
dependence of $\xi$
on the underlying point configuration~$\mu$. The difficulty comes from
the fact that adding $\d_{x_1}+\cdots+\d_{x_k}$ to the Poisson process
$\eta$ amounts not only to adding $\d_{f(x_1,\ldots,x_k)}$ to
$\xi(\eta)$ but also all atoms of the form $f(x_{i_1},\ldots,
x_{i_\ell},\tilde{x}_{\ell+1}, \ldots, \tilde{x}_k)$ with
$\ell\in\{1,\ldots,k\}$, pairwise different indices $i_1,\ldots
,i_\ell\in\{1,\ldots,k\}$ and $(\tilde{x}_{\ell+1}, \ldots,
\tilde{x}_k)\in\eta_{\neq}^{k-\ell}$. We denote by
$\hat{\xi}(x_1,\ldots,x_k,\eta)$ the collection of these extra atoms.
The difference in \eqref{eqPPPconv5} is now decomposed as
%
\begin{eqnarray}
\label{eqDecompositionIdea}
&& \bE\bigl[ \bigl( F\bigl(\xi(\eta)+ \hat
{\xi}(x_1,
\ldots,x_k,\eta)\bigr)-F\bigl(\xi(\eta)\bigr) \bigr)\nonumber
\\
&&\hspace*{-9pt}\qquad{} + \bigl(F\bigl(
\xi(\eta)\bigr)-F\bigl(\xi(\eta)+\delta_{f(x_1,\ldots,x_k)}\bigr) \bigr)
\\
&&\hspace*{-9pt}\qquad{} + \bigl(F\bigl(\xi(\eta)+\delta_{f(x_1,\ldots,x_k)}\bigr)-F\bigl
(\xi(\eta)+
\hat{\xi}(x_1,\ldots,x_k,\eta)+\d_{f(x_1,\ldots,x_k)}\bigr)
\bigr) \bigr] .\hspace*{-20pt}\nonumber
\end{eqnarray}
The middle term in \eqref{eqDecompositionIdea} contributes to the
total variation distance of the intensity measures in \eqref
{eqPPPconv6} in Theorem
\ref{thmMain}. Since $F$ is Lipschitz, the expectation and the
integral with respect to $x_1,\ldots,x_k$ of the first and
the third term in \eqref{eqDecompositionIdea} are bounded (up to a
constant) by
\[
\bE\int_{\XX^k} \hat{\xi}(x_1,
\ldots,x_k,\eta) (\YY) \bK^k\bigl(\dint(x_1,
\ldots,x_k)\bigr) ,
\]
which in turn is bounded by $\bE\xi(\YY)^2-\bE\xi(\YY)-(\bE\xi
(\YY))^2$ and $r(\dom f)$. This effect contributes to the second term
of the bounds in Theorem~\ref{thmMain}. For $k=1$, only the middle
term in \eqref{eqDecompositionIdea} is present. This explains why,
for $k=1$, the Kantorovich--Rubinstein distance between the
transformation of a Poisson process (which is again a Poisson process)
and a second Poisson process is bounded by the total variation distance
of the intensity measures, and the second term in \eqref{eqPPPconv6}
in Theorem~\ref{thmMain} vanishes.

Throughout this section we use the same notation as in Section~\ref
{subsecGeneralEstimate}. Moreover, let $[k]$ be shorthand for $\{
1,\ldots,k\}$. For\vspace*{1pt} $x=(x_1,\ldots,x_k)\in\XX^k$, $I=\{i_1,\ldots
,i_{\llvert I\rrvert}\}\subset[k]$ and $z=(z_1,\ldots,z_{k-\llvert
I\rrvert})\in\XX^{k-\llvert I\rrvert}$,
let $(x_I,z)=(x_{i_1},\ldots,x_{i_{\llvert I\rrvert}},\break z_1,\ldots
,z_{k-\llvert I\rrvert})$. We
prepare the proof of Theorem~\ref{thmMain} with the following lemma.

%
\begin{lemma}\label{lemProductFormula}
Let the assumptions of Theorem~\ref{thmMain} prevail. If $\xi$ is
induced by a Poisson process, then
\begin{eqnarray*}
\bE\xi(\YY)^2 &=& \frac{1}{k!} \sum_{I\subset[k]}
\frac
{1}{(k-\llvert I\rrvert)!} \int_{\XX^k} \int_{\XX^{k-\llvert
I\rrvert}}
 {\mathbf1}\bigl((x_1,\ldots,x_k)\in\dom f\bigr)
\\
&&\hspace*{126pt}{} \times{\mathbf1}\bigl((x_I,z)\in\dom f\bigr) \bK^{k-\llvert
I\rrvert}(
\dint z)
\\
&&\hspace*{126pt}{} \times\bK^k\bigl(\dint(x_1,\ldots,x_k)
\bigr) .
\end{eqnarray*}
If $\xi$ is derived from a binomial process of $n$ points, then
\begin{eqnarray*}
\bE\xi(\YY)^2 &=& \frac{1}{k!} \sum_{I\subset[k]}
\frac
{(n)_{2k-\llvert I\rrvert}}{(k-\llvert I\rrvert)!} \int_{\XX^k}
\int_{\XX^{k-\llvert I\rrvert}}
{\mathbf1}\bigl((x_1,\ldots,x_k)\in\dom f\bigr)
\\
&&\hspace*{126pt}{} \times{\mathbf1}\bigl((x_I,z)\in\dom f\bigr) \bK_1^{k-\llvert I
\rrvert}(
\dint z)
\\
&&\hspace*{126pt}{} \times\bK_1^k\bigl(\dint(x_1,
\ldots,x_k)\bigr) .
\end{eqnarray*}
\end{lemma}
\begin{pf}
We have that
\begin{eqnarray*}
\xi(\YY)^2 & =&\frac{1}{(k!)^2} \biggl( \sum
_{(x_1,\ldots,x_k)\in
\mu^k_{\neq}} {\mathbf1}\bigl( (x_1,\ldots,x_k)
\in\dom f \bigr) \biggr)^2
\\
& =& \frac{1}{(k!)^2}\sum_{I\subset[k]} \sum
_{(x_1,\ldots,x_k,z)\in
\mu^{2k-\llvert I \rrvert}_{\neq}} \frac{k!}{(k-\llvert I \rrvert
)!} {\mathbf1}\bigl((x_1,\ldots,x_k)\in\dom f\bigr)
\\
&&\hspace*{163pt}{}\times {\mathbf1}\bigl((x_I,z)\in\dom f\bigr),
\end{eqnarray*}
where we have used that two points occurring in different sums can be
either equal or distinct and that $\dom f$ is symmetric. Now the
multivariate Mecke \eqref{eqMecke} and its binomial analogue \eqref
{eqMomentMeasureBinomial} complete the proof.
\end{pf}

\begin{pf*}{Proof of Theorem~\ref{thmMain}}
Throughout this proof we write $\xi(\eta)$ and $\xi(\beta_n)$ to
emphasize the dependence of $\xi$ on the underlying point process. Whenever
we do not need special properties of $\eta$ or $\beta_n$, we write
$\xi(\mu)$ with the dummy variable $\mu$ standing for either $\eta$
or $\beta_n$. As discussed in Remark~\ref{remMain}(ii), we can
assume for the Poisson case that $\bL(\YY)<\infty$ and hence that
$\xi(\eta)$ is almost surely finite since \eqref{eqPPPconv6} is
obviously true otherwise. For an underlying binomial process it is
sufficient to consider only the case $n\geq k$ since, otherwise, the
statement is obviously true as explained in Remark~\ref{remMain}(i).

Lemma~\ref{lemPPPconv1} says that for $h\in\cC$ and $\o\in
\widehat{\N}_\YY$,
%
\begin{equation}
\label{eqSteinIdentity} \bE h(\zeta)- h(\o)=\int_0^\infty
LP_sh(\o) \,\dint s .
\end{equation}
The Stein-type identity
\eqref{eqSteinIdentity} is the starting point for our proof.
Combining \eqref{eqSteinIdentity}
with the representation of the generator $L$ in
\eqref{eqGeneratorL}, choosing $\omega=\xi(\mu)$ and taking
expectations results in the following:
%
\begin{eqnarray}
\label{eqdifferenceExpectations} \bE h(\zeta) - \bE h\bigl(\xi(\mu
)\bigr) &=& \bE\int
_0^\infty LP_sh\bigl(\xi(\mu)\bigr)\,
\dint s\nonumber
\\
&=&\bE\int_0^\infty\!\!\int_{\YY}
\bigl(P_sh\bigl(\xi(\mu)+\d_y\bigr)-P_sh
\bigl(\xi(\mu)\bigr) \bigr) \bM(\dint y) \,\dint s
\\
&&{} +\bE\int_0^\infty\sum
_{y\in\xi(\mu)} \bigl(P_sh\bigl(\xi(\mu)-
\d_y\bigr)-P_sh\bigl(\xi(\mu)\bigr) \bigr) \,\dint s .\nonumber
\end{eqnarray}
Let us denote the first and the second term on the right-hand side
by $T_{1,\mu}$ and $T_{2,\mu}$, respectively. By Fubini's theorem and
the definition of $\xi(\mu)$,
we obtain that
\[
T_{2,\mu} =\frac{1}{k!}\int_0^\infty
\bE\sum_{(x_1,\ldots,x_k)\in\mu_{\neq
}^k\cap
\dom f} \bigl(P_sh\bigl(\xi(\mu)-
\d_{f(x_1,\ldots,x_k)}\bigr)-P_sh\bigl(\xi(\mu)\bigr) \bigr) \,\dint s .
\]
By the multivariate Mecke formula \eqref{eqMecke} and its
analogue \eqref{eqMomentMeasureBinomial} for binomial processes, we
see that
\begin{eqnarray*}
T_{2,\eta} &=& \frac{1}{k!}\int_0^\infty\!\!\int_{\dom f}\bE\bigl[ P_sh\bigl(\xi(\eta+
\d_{x_1}+\cdots+\d_{x_k})-\d_{f(x_1,\ldots,x_k)}\bigr)
\\[-1.5pt]
&&\hspace*{47pt}\hspace*{69pt}{}-P_sh\bigl(\xi(\eta+\d_{x_1}+\cdots+
\d_{x_k})\bigr) \bigr]
\\[-1.5pt]
&&\hspace*{56pt}{}\times \bK^k\bigl(\dint(x_1,
\ldots,x_k)\bigr) \,\dint s
\end{eqnarray*}
and
\begin{eqnarray*}
T_{2,\beta_n} &=& \frac{(n)_k}{k!}\int_0^\infty\!\!\int_{\dom f}\bE\bigl[ P_sh\bigl(\xi(
\beta_{n-k}+\d_{x_1}+\cdots+\d_{x_k})-\d
_{f(x_1,\ldots,x_k)}\bigr)
\\[-1.5pt]
&&\hspace*{127pt}{} -P_sh\bigl(\xi(\beta_{n-k}+
\d_{x_1}+\cdots+\d_{x_k})\bigr) \bigr]
\\[-1.5pt]
&&\hspace*{66pt}{}\times  \bK^k_1
\bigl(\dint(x_1,\ldots,x_k)\bigr) \,\dint s .
\end{eqnarray*}
Let us write $\hat\xi(x_1,\ldots,x_k,\mu)$ for the point process
\[
\hat\xi(x_1,\ldots,x_k,\mu):=\sum
_{\varnothing\neq I\subsetneq
[k], z\in\mu_{\neq}^{k-\llvert I \rrvert}} \frac{1}{(k-\llvert I
\rrvert)!} {\mathbf1}\bigl((x_I,z)\in
\dom f\bigr) \d_{f(x_I,z)} %
\]
on $\YY$, where $\subsetneq$
denotes proper set-inclusion and where the notation $(x_I,z)$ has been
introduced before Lemma~\ref{lemProductFormula} above. Then
\begin{eqnarray*}
T_{2,\eta} &=&\frac{1}{k!}\int_0^\infty\!\!\int_{\dom f}\bE\bigl[ P_sh\bigl(\xi(\eta)+\hat
\xi(x_1,\ldots,x_k,\eta)\bigr)
\\[-1.5pt]
&&\hspace*{69pt}{} -P_sh\bigl(\xi(\eta)+\hat\xi(x_1,
\ldots,x_k,\eta)+\d_{f(x_1,\ldots,x_k)}\bigr) \bigr]
\\[-1.5pt]
&&\hspace*{55pt}{}\times  \bK^k
\bigl(\dint(x_1,\ldots,x_k)\bigr)\,\dint s
\\[-1.5pt]
&=&-\frac{1}{k!}\int_0^\infty\!\!\int
_{\dom
f}\bE\bigl[P_sh\bigl(\xi(\eta)+
\d_{f(x_1,\ldots,x_k)}\bigr)-P_sh\bigl(\xi(\eta)\bigr) \bigr]
\\[-1.5pt]
&&\hspace*{64pt}{}\times \bK^k\bigl(\dint(x_1,\ldots,x_k)\bigr)\,\dint s
\\[-1.5pt]
&&{} +\frac{1}{k!}\int_0^\infty\!\!\int
_{\dom f}\bE\bigl[P_sh\bigl(\xi(\eta)+\hat
\xi(x_1,\ldots,x_k,\eta)\bigr)
\\[-1.5pt]
&&\hspace*{81pt}{} -P_sh\bigl(
\xi(\eta)\bigr)+P_sh\bigl(\xi(\eta)+\d_{f(x_1,\ldots,x_k)}\bigr)
\\[-1.5pt]
&&\hspace*{81pt}{} -P_sh\bigl(\xi(\eta)+\hat\xi(x_1,
\ldots,x_k,\eta)+\d_{f(x_1,\ldots
,x_k)}\bigr) \bigr]
\\[-1.5pt]
&&\hspace*{68pt}{}\times \bK^k
\bigl(\dint(x_1,\ldots,x_k)\bigr)\,\dint s
\\[-1.5pt]
&=:&\hat T_{2,\eta} + R_\eta
\end{eqnarray*}
and
\begin{eqnarray*}
T_{2,\beta_n} &=&\frac{(n)_k}{k!}\int_0^\infty\!\!\int_{\dom f}\bE\bigl[ P_sh\bigl(\xi(
\beta_{n-k})+\hat\xi(x_1,\ldots,x_k,
\beta_{n-k})\bigr)
\\
&&\hspace*{79pt}{} -P_sh\bigl(\xi(\beta_{n-k})+
\hat\xi(x_1,\ldots,x_k,\beta_{n-k})+
\d_{f(x_1,\ldots,x_k)}\bigr) \bigr]
\\
&&\hspace*{66pt}{}\times  \bK_1^k\bigl(
\dint(x_1,\ldots,x_k)\bigr)\,\dint s
\\
& =&-\frac{(n)_k}{k!}\int_0^\infty\!\!\int
_{\dom
f}\bE\bigl[P_sh\bigl(\xi(
\beta_{n-k})+\d_{f(x_1,\ldots,x_k)}\bigr)-P_sh\bigl(\xi(
\beta_{n-k})\bigr) \bigr]
\\
&&\hspace*{75pt}{}\times  \bK_1^k\bigl(
\dint(x_1,\ldots,x_k)\bigr)\,\dint s
\\
&&{}+\frac{(n)_k}{k!}\int_0^\infty\!\!\int
_{\dom f}\bE\bigl[P_sh\bigl(\xi(
\beta_{n-k})+\hat\xi(x_1,\ldots,x_k,
\beta_{n-k})\bigr)
\\
&&\hspace*{92pt}{} -P_sh\bigl(\xi(\beta_{n-k})
\bigr)
+P_sh\bigl(\xi(\beta_{n-k})+\d
_{f(x_1,\ldots,x_k)}\bigr)
\\
&&\hspace*{92pt}{} -P_sh\bigl(\xi(\beta_{n-k})+
\hat\xi(x_1,\ldots,x_k,\beta_{n-k})
\\
&&\hspace*{144pt}\hspace*{92pt}{} + \d_{f(x_1,\ldots,x_k)}\bigr) \bigr]
\\
&&\hspace*{79pt}{}\times  \bK_1^k\bigl(
\dint(x_1,\ldots,x_k)\bigr)\,\dint s
\\
&=:&\hat T_{2,\beta_n} + R_{\beta_n} .
\end{eqnarray*}
Together with \eqref{eqdifferenceExpectations} and the formulas for
$\bL$ in \eqref{eqLPoisson} and \eqref{eqLbinomial}, we see that
\[
\bE h(\zeta) - \bE h\bigl(\xi(\eta)\bigr) = \int_0^\infty\!\!\int_{\YY}\bE\bigl[D_yP_sh\bigl(
\xi(\eta)\bigr) \bigr] (\bM-\bL) (\dint y)\,\dint s +R_\eta%
\]
and
\begin{eqnarray*}
&& \bE h(\zeta) - \bE h\bigl(\xi(\beta_n)\bigr)
\\
&&\qquad = \int _0^\infty\!\!\int_{\YY}\bE
\bigl[D_yP_sh\bigl(\xi(\beta_n)\bigr)
\bigr] (\bM-\bL) (\dint y)\,\dint s
\\
&&\quad\qquad{} + \int_0^\infty\!\!\int_{\YY}\bE
\bigl[D_yP_sh\bigl(\xi(\beta_n)\bigr)
\bigr]-\bE\bigl[D_yP_sh\bigl(\xi(\beta_{n-k})
\bigr) \bigr] \bL(\dint y)\,\dint s +R_{\beta_n} .
\end{eqnarray*}
We now determine the remainder terms $R_\eta$ and
$R_{\beta_n}$. For $(x_1,\ldots,x_k)\in\dom f$ let us define
$\widetilde{h}_{x_1,\ldots,x_k}\dvtx \widehat{\N}_\YY\to\mathbb{R}$ by
\[
\widetilde{h}_{x_1,\ldots,x_k}(\mu)=\tfrac{1}{2} \bigl(h(\mu)-h(\mu+
\delta_{f(x_1,\ldots,x_k)}) \bigr) .
\]
We can then rewrite $R_\eta$ and $R_{\beta_n}$ as
\begin{eqnarray*}
R_\eta&=&\frac{2}{k!} \int_0^\infty\!\!\int_{\dom f} \bE\bigl[P_s\widetilde{h}_{x_1,\ldots,x_k}
\bigl(\xi(\eta)+\hat\xi(x_1,\ldots,x_k,\eta)\bigr)
-P_s\widetilde{h}_{x_1,\ldots,x_k}\bigl(\xi(\eta)\bigr) \bigr]
\\
&&\hspace*{55pt}{}\times
\bK^k\bigl(\dint(x_1,\ldots,x_k)\bigr)\,\dint s
\end{eqnarray*}
and
\begin{eqnarray*}
R_{\beta_n} &=& \frac{2(n)_k}{k!} \int_0^\infty\!\!\int_{\dom f} \bE\bigl[P_s\widetilde{h}_{x_1,\ldots,x_k}
\bigl(\xi(\beta_{n-k}) +\hat\xi(x_1,\ldots,x_k,
\beta_{n-k})\bigr)
\\[1.5pt]
&&\hspace*{172pt}{}-P_s
\widetilde{h}_{x_1,\ldots,x_k}\bigl(\xi(\beta_{n-k})\bigr) \bigr]
\\[1.5pt]
&&\hspace*{71pt}{}\times
\bK_1^k\bigl(\dint(x_1,\ldots,x_k)
\bigr)\,\dint s .
\end{eqnarray*}
Because of $\widetilde{h}_{x_1,\ldots,x_k}\in\cC$, we obtain by the
definition of the semi-group $(P_s)_{s\geq0}$ in \eqref
{eqSemigroupPt} and Lemma
\ref{lemPreparationForOtherLemmas} that
\begin{eqnarray*}
\llvert R_\eta\rrvert& \leq&\frac{2}{k!} \int _0^\infty\!\!\int_{\dom f}
e^{-s} \bE\hat{\xi}(x_1,\ldots,x_k,\eta) (\YY)
\bK^k\bigl(\dint(x_1,\ldots,x_k)\bigr)\,\dint s
\\[1.5pt]
& =& \frac{2}{k!} \int_{\dom f} \bE\hat{
\xi}(x_1,\ldots,x_k,\eta) (\YY) \bK^k\bigl(
\dint(x_1,\ldots,x_k)\bigr)
\end{eqnarray*}
and
\begin{eqnarray*}
\llvert R_{\beta_n}\rrvert& \leq&\frac{2(n)_k}{k!} \int _0^\infty\!\!\int_{\dom f}
e^{-s} \bE\hat{\xi}(x_1,\ldots,x_k,
\beta_{n-k}) (\YY) \bK_1^k\bigl(
\dint(x_1,\ldots,x_k)\bigr)\,\dint s
\\[1.5pt]
& \leq&\frac{2(n)_k}{k!} \int_{\dom f} \bE\hat{
\xi}(x_1,\ldots,x_k,\b_{n-k}) (\YY)
\bK^k_1\bigl(\dint(x_1,\ldots,x_k)
\bigr) .
\end{eqnarray*}
Now, from the Mecke formula \eqref{eqMecke} and its analogue \eqref
{eqMomentMeasureBinomial} for binomial processes, it
follows that
\begin{eqnarray*}
&& \bE\hat\xi(x_1,\ldots,x_k,\eta) (\YY)
\\[1.5pt]
&&\qquad = \bE\sum
_{\varnothing\neq
I\subsetneq[k], z\in\eta_{\neq}^{k-\llvert I \rrvert}} \frac
{1}{(k-\llvert I \rrvert)!} {\mathbf1}\bigl(f(x_I,z)
\in\dom f\bigr)
\\[1.5pt]
&&\qquad = \sum_{\varnothing\neq I\subsetneq[k]} \frac{1}{(k-\llvert I \rrvert
)!} \int
_{\XX^{k-\llvert I \rrvert}} {\mathbf1}\bigl((x_I,z)\in\dom f\bigr)
\bK^{k-\llvert I \rrvert}(\dint z) 
\end{eqnarray*}
and
%
\begin{eqnarray}
\label{eqxiHatBinomial}
&& \bE\hat\xi(x_1,\ldots,x_k,
\beta_{n-k}) (\YY)\nonumber
\\[1.5pt]
&&\qquad = \bE\sum_{\varnothing\neq I\subsetneq[k], z\in\beta
_{n-k,\neq}^{k-\llvert I \rrvert}}
\frac{1}{(k-\llvert I \rrvert)!} {\mathbf1}\bigl(f(x_I,z)\in\dom
f\bigr)
\\[1.5pt]
&&\qquad = \sum_{\varnothing\neq I\subsetneq[k]} \frac{(n-k)_{k-\llvert I
\rrvert}}{(k-\llvert I \rrvert)!} \int
_{\XX^{k-\llvert I \rrvert}} {\mathbf1}\bigl((x_I,z)\in\dom f\bigr)
\bK_1^{k-\llvert I \rrvert}(\dint z) .\nonumber
\end{eqnarray}
Together with Lemma~\ref{lemProductFormula}, we obtain
%
\begin{eqnarray}
\label{eqbound1Poisson} \llvert R_\eta\rrvert& \leq&
\frac{2}{k!}\int_{\XX^k }\sum_{\varnothing\neq I\subsetneq[k]}
\frac{1}{(k-\llvert I \rrvert)!} \int_{\XX
^{k-\llvert I \rrvert}}{\mathbf1}\bigl((x_1,
\ldots,x_k)\in\dom f\bigr)\nonumber
\\
&&\hspace*{137pt}{} \times{\mathbf1}\bigl((x_I,z)\in\dom f\bigr)
\bK^{k-\llvert I \rrvert}(\dint z)
\nonumber\\[-8pt]\\[-8pt]\nonumber
&&\hspace*{137pt}{} \times\bK^k\bigl(\dint(x_1,
\ldots,x_k)\bigr)
\\
&=& 2 \bigl(\bE\xi(\YY)^2-\bL(\YY)-\bL(\YY)^2 \bigr)= 2
\bigl(\bE\xi(\YY)^2-\bE\xi(\YY)-\bigl(\bE\xi(\YY)\bigr)^2
\bigr)\nonumber
\end{eqnarray}
and
%
\begin{eqnarray}\label{eqbound1Binomial}
\llvert R_{\beta_n}\rrvert& \leq&\frac
{2}{k!}\int
_{\XX^k }\sum_{\varnothing\neq
I\subsetneq[k]}
\frac{(n)_k (n-k)_{k-\llvert I \rrvert}}{(k-\llvert I \rrvert)!}
\int_{\XX
^{k-\llvert I \rrvert}}{\mathbf1}\bigl((x_1,
\ldots,x_k)\in\dom f\bigr)\nonumber
\\
&&\hspace*{167pt}{}\times{\mathbf1}\bigl((x_I,z)\in\dom f\bigr)
\bK^{k-\llvert I \rrvert}_1(\dint z)\nonumber
\\
&&\hspace*{167pt}{}\times  \bK^k_1\bigl(
\dint(x_1,\ldots,x_k)\bigr)
\\
&=& 2 \biggl(\bE\xi(\YY)^2-\bL(\YY)-\frac{(n-k)_{k}}{(n)_k}\bL(\YY
)^2 \biggr)\nonumber
\\
& =& 2 \biggl(\bE\xi(\YY)^2-\bE\xi(\YY)-
\frac
{(n-k)_{k}}{(n)_k}\bigl(\bE\xi(\YY)\bigr)^2 \biggr) .\nonumber
\end{eqnarray}
The inequalities in \eqref{eqbound1Poisson} and \eqref
{eqbound1Binomial} together with the definition of $r(\dom f)$ imply that
%
\begin{equation}
\label{eqbound1PoissonBinomial} \llvert R_\eta\rrvert\leq\frac
{2^{k+1}}{k!} r(\dom f)
\quad\mbox{and}\quad\llvert R_{\beta_n}\rrvert\leq\frac{2^{k+1}}{k!}
r(\dom
f) .
\end{equation}
Next, it follows from Lemma~\ref{lemPreparationForOtherLemmas} that
for $s\geq0$,
%
\begin{equation}
\label{eqEDPh} \bigl\llvert\bE D_y P_s h\bigl(\xi(
\mu)\bigr) \bigr\rrvert\leq\bE\bigl[ \bigl\llvert P_sh\bigl(\xi(
\mu)+\delta_y\bigr) - P_sh\bigl(\xi(\mu)\bigr)\bigr
\rrvert\bigr]\leq e^{-s} .
\end{equation}
For $y_1,y_2\in\YY$ and $\tilde{\xi}\in\widetilde{\mathrm{N}}_{\YY
}$ we have $\dTV(\tilde{\xi}+\delta_{y_1},\tilde{\xi}+\delta
_{y_2})\leq1$ so that $h\in\cC$ leads to
\[
\bigl\llvert D_{y_1}h(\tilde{\xi})-D_{y_2}h(\tilde{\xi})
\bigr\rrvert=\bigl\llvert h(\tilde{\xi}+\delta_{y_1})-h(\tilde{\xi}+
\delta_{y_2})\bigr\rrvert\leq1 . %
\]
Together with Lemma~\ref{lemDPt}, we obtain that
%
\begin{eqnarray}
\label{eqDiffEDPh}
&& \bigl\llvert\bE D_{y_1}P_sh\bigl(\xi(\mu)
\bigr)-\bE D_{y_2}P_sh\bigl(\xi(\mu)\bigr)\bigr\rrvert
\nonumber\\[-8pt]\\[-8pt]\nonumber
&&\qquad =e^{-s} \bigl\llvert\bE P_s(D_{y_1}h-D_{y_2}h)
\bigl(\xi(\mu)\bigr)\bigr\rrvert\leq e^{-s}
\end{eqnarray}
for all $y_1,y_2\in\YY$ and $s\geq0$. The estimates in \eqref
{eqEDPh} and \eqref{eqDiffEDPh} show that
%
\begin{eqnarray}\label{eqbound2}
\biggl\llvert\int_0^\infty\!\!\int
_\YY\bE\bigl[D_yP_sh\bigl(\xi(
\mu)\bigr) \bigr] (\bM-\bL) (\dint y)\,\dint s \biggr\rrvert
&\leq&\dTV(\bM,\bL)\int_0^\infty e^{-s}\,
\dint s
\nonumber\\[-8pt]\\[-8pt]\nonumber
& \leq&\dTV(\bM,\bL) .
\end{eqnarray}
Combining \eqref{eqbound1Poisson} and \eqref
{eqbound1PoissonBinomial} with \eqref{eqbound2} completes the
proof of the Poisson case.

When considering a binomial process, we
additionally need to take care of the term
\[
\int_0^\infty\!\!\int_{\YY}\bE
\bigl[D_yP_sh\bigl(\xi(\beta_n)\bigr)
\bigr]-\bE\bigl[D_yP_sh\bigl(\xi(\beta_{n-k})
\bigr) \bigr] \bL(\dint y)\,\dint s . %
\]
For this, we use Lemma~\ref{lemDPt}, the fact that
$\frac{1}{2}D_yh\in\cC$ whenever $h\in\cC$ and Lemma~\ref
{lemPreparationForOtherLemmas} to obtain that
\begin{eqnarray*}
&&\bigl\llvert\bE\bigl[D_yP_sh\bigl(\xi(
\beta_n)\bigr) \bigr]-\bE\bigl[D_yP_sh
\bigl(\xi(\beta_{n-k})\bigr) \bigr] \bigr\rrvert
\\
&&\qquad  \leq\int_{\XX^k} \bigl\llvert\bE\bigl[D_yP_sh
\bigl(\xi(\beta_{n-k}+\delta_{x_1}+\cdots+
\delta_{x_k})\bigr) \bigr]-\bE\bigl[D_yP_sh
\bigl(\xi(\beta_{n-k})\bigr) \bigr] \bigr\rrvert
\\
&&\hspace*{47pt}{} \times \bK_1^k
\bigl(\dint(x_1,\ldots,x_k)\bigr)
\\
&&\qquad  = \int_{\XX^k} e^{-s} \bigl\llvert\bE
\bigl[P_s(D_yh) \bigl(\xi(\beta_{n-k})+\hat{
\xi}(x_1,\ldots,x_k,\beta_{n-k})+
\delta_{f(x_1,\ldots,x_k)}\bigr)
\\
&&\hspace*{225pt}{} -P_s(D_yh) \bigl(\xi(\beta
_{n-k})\bigr) \bigr] \bigr\rrvert
\\
&&\hspace*{46pt}{}\times  \bK_1^k\bigl(\dint(x_1,
\ldots,x_k)\bigr)
\\
&&\qquad \leq\frac{1}{n^k} \int_{\XX^k} 2e^{-2s} \bigl(
\bE\hat{\xi}(x_1,\ldots,x_k,\b_{n-k}) (\YY)+{
\mathbf1}\bigl((x_1,\ldots,x_k)\in\dom f\bigr) \bigr)
\\
&&\hspace*{61pt}{}\times \bK^k\bigl(\dint(x_1,\ldots,x_k)\bigr)
\end{eqnarray*}
for any $s\geq0$. It follows from \eqref{eqxiHatBinomial} and
$(n-k)_{k-\llvert I \rrvert}\leq n^{k-\llvert I \rrvert}$ that
\begin{eqnarray*}
&& \frac{1}{n^k} \int_{\XX^k} \bE\hat{
\xi}(x_1,\ldots,x_k,\b_{n-k}) (\YY)
\bK^k\bigl(\dint(x_1,\ldots,x_k)\bigr)
\\
&&\qquad  \leq\frac{1}{n^k} \int_{\XX^k} \sum
_{\varnothing\neq I\subsetneq
[k]} \frac{1}{(k-\llvert I \rrvert)!} \int_{\XX^{k-\llvert I
\rrvert}} {
\mathbf1}\bigl((x_I,z)\in\dom f\bigr) \bK^{k-\llvert I \rrvert
}(\dint z)
\\
&&\hspace*{172pt}{}\times \bK^k\bigl(\dint(x_1,\ldots,x_k)\bigr)
\\
&&\qquad  = \frac{1}{n^k} \sum_{\varnothing\neq I\subsetneq[k]}
\frac
{1}{(k-\llvert I \rrvert)!} \int_{\XX^k} {\mathbf1}\bigl((x_1,
\ldots,x_k)\in\dom f\bigr)
\\
&&\hspace*{143pt}{}\times  \bK^{k}\bigl(
\dint(x_1,\ldots,x_k)\bigr)\bK(\XX)^{\llvert k \rrvert-\llvert I
\rrvert}
\\
&&\qquad  \leq\frac{(2^{k}-2)}{n} \int_{\XX^k} {\mathbf1}
\bigl((x_1,\ldots,x_k)\in\dom f\bigr) \bK^{k}
\bigl(\dint(x_1,\ldots,x_k)\bigr) .
\end{eqnarray*}
Now, \eqref{eqLbinomial} implies that
\[
\int_{\XX^k} {\mathbf1}\bigl((x_1,
\ldots,x_k)\in\dom f\bigr) \bK^k\bigl(\dint
(x_1,\ldots,x_k)\bigr)=\frac{k! n^k}{(n)_k} \bL(\YY) \leq
k! e^k \bE\xi(\YY) , %
\]
where we have used that $n^k/(n)_k\leq k^k/k!\leq e^k$ for $n\geq k$.
Hence, using that $2^k e^k\leq6^k$, we find
\begin{eqnarray*}
\int_0^\infty\!\!\int_{\YY} \bigl
\llvert\bE\bigl[D_yP_sh\bigl(\xi(\beta_n)
\bigr) \bigr]-\bE\bigl[D_yP_sh\bigl(\xi(
\beta_{n-k})\bigr) \bigr]\bigr\rrvert\bL(\dint y)\,\dint s
&\leq& 6^k k! \frac{\bL(\YY)^2}{n}
\\
&=& 6^k k! \frac
{(\bE\xi(\YY))^2}{n} .
\end{eqnarray*}
Together with \eqref{eqbound1Binomial}, \eqref
{eqbound1PoissonBinomial} and \eqref{eqbound2} this completes the
proof in the binomial case.
\end{pf*}

%
\begin{remark}
Bounds for the total variation distance between $\xi$ and $\zeta$
that are similar to the bounds for the Kantorovich--Rubinstein distance
in Theorem~\ref{thmMain} can be deduced from Theorem 2.6 in \cite
{BarbourBrown}. This result implies that
\begin{eqnarray*}
\dTVRandom\bigl(\xi(\eta),\zeta\bigr) &\leq& 2\dTV(\bL,\bM)
\\
&&{}+\frac{2}{k!}\int_{\dom f} \bE\dTV\bigl(\xi(\eta),\xi(
\eta+\delta_{x_1}+\cdots+\delta_{x_k})-\delta_{f(x_1,\ldots,x_k)}
\bigr)
\\
&&\hspace*{51pt}{}\times  \bK^k\bigl(\dint(x_1,\ldots,x_k)\bigr)
\end{eqnarray*}
and
\begin{eqnarray*}
&& \dTVRandom\bigl(\xi(\beta_n) ,\zeta\bigr)
\\
&&\qquad \leq 2\dTV(\bL,\bM)
\\
&&\quad\qquad{}+\frac
{2(n)_k}{k!}\int_{\dom f} \bE\dTV\bigl(\xi(
\beta_n),\xi(\beta_{n-k}+\delta_{x_1}+\cdots+
\delta_{x_k})-\delta_{f(x_1,\ldots
,x_k)}\bigr)
\\
&&\hspace*{99pt}{}\times \bK_1^k\bigl(\dint(x_1,\ldots,x_k)\bigr) .
\end{eqnarray*}
Since the integrands are bounded by
\[
\bE\hat{\xi}(x_1,\ldots,x_k,\eta) (\YY)\quad\mbox{and}
\quad\bE\hat{\xi}(x_1,\ldots,x_k,\beta_{n-k}) (
\YY) + \bE\dTV\bigl(\xi(\beta_n),\xi(\beta_{n-k})\bigr) ,
\]
respectively, the integrals on the right-hand sides can be controlled
as in the proof of Theorem~\ref{thmMain} above.
\end{remark}

\section{Applications}\label{secApplications}

\subsection{Poisson approximation of U-statistics}\label{secUstatistics}

In this subsection we present a first application of Theorem~\ref
{thmMain} to
U-statistics of Poisson or binomial processes. Let $(\XX,\cX)$ and
$(\YY,\cY)$ be two lcscH spaces, and let for some fixed integer
$k\geq1$, $f_t\dvtx\XX^k\to\YY$, $t\geq1$, be symmetric measurable
functions. Furthermore, for
a $\sigma$-finite measure $\bK$ and a probability measure $\bK_1$ on
$\XX$, we denote by $\eta_t$ a Poisson process with intensity measure
$\bK_t:=t\bK$, $t\geq1$, and by $\beta_t$, $t\geq1$, a binomial
process of $\lceil t\rceil$ points with
intensity measure $\bK_t:=\lceil t\rceil\bK_1$, respectively. If
$\mu_t$ is either $\eta_t$ or $\beta_t$ and if $B$ is a measurable
subset of $\YY$, we define the U-statistics
\[
S_t(B):=\frac{1}{k!}\sum_{(x_1,\ldots,x_k)\in\mu_{t,\neq}^k}{\mathbf
1}\bigl(f_t(x_1,\ldots,x_k)\in B\bigr) ,
\qquad t\geq1 , %
\]
which count the number of $k$-tuples $(x_1,\ldots,x_k)\in\mu_{t,\neq
}^k$ for which
$f_t(x_1,\ldots,\break x_k)\in B$. To compare $S_t(B)$ with a Poisson random
variable, we define
\begin{eqnarray*}
r_t(B)&:=&\max_{1\leq\ell\leq k-1} \int_{\XX^\ell}
\biggl(\int_{\XX^{k-\ell}} {\mathbf1}\bigl(f_t(x_1,
\ldots,x_k)\in B\bigr) \bK_t^{k-\ell
}\bigl(
\dint(x_{\ell+1},\ldots,x_k)\bigr) \biggr)^2
\\
&&\hspace*{51pt}{}\times \bK_t^{\ell}\bigl(\dint(x_1,\ldots,x_\ell)
\bigr) %
\end{eqnarray*}
if $k>1$ and $r_t(B):=0$ if $k=1$.

%
\begin{theorem}\label{thmUstatistics}
Let $B\in\cY$, and let $Z$ be a Poisson distributed random variable
with mean $\lambda\in[0,\infty)$. Suppose that $\bE S_t(B)^2<\infty
$. If $S_t(B)$ is induced
by a Poisson process $\eta_t$ with $t\geq1$, then
\begin{eqnarray*}
\dW\bigl(S_t(B),Z\bigr) & \leq&\bigl\llvert\bE S_t(B)-
\lambda\bigr\rrvert+2 \bigl(\bE S_t(B)^2-\bE
S_t(B)-\bigl(\bE S_t(B)\bigr)^2 \bigr)
\\
& \leq&\bigl\llvert\bE S_t(B)-\lambda\bigr\rrvert+
\frac{2^{k+1}}{k!} r_t(B) .
\end{eqnarray*}
If $S_t(B)$ is induced by a binomial process $\beta_t$ with $t\geq1$, then
\begin{eqnarray*}
&&\dW\bigl(S_t(B),Z\bigr) \\
&&\qquad \leq\bigl\llvert\bE S_t(B)-
\lambda\bigr\rrvert+2 \biggl(\bE S_t(B)^2-\bE
S_t(B)-\frac{(\lceil t\rceil-k)_k}{(\lceil t\rceil)_k}\bigl(\bE S_t(B)
\bigr)^2 \biggr)
\\
&&\qquad\quad{} +\frac{6^k k!}{t} \bigl(\bE S_t(B)\bigr)^2
\\
& &\qquad\leq\bigl\llvert\bE S_t(B)-\lambda\bigr\rrvert+\frac{2^{k+1}}{k!}
r_t(B)+\frac{6^k
k!}{t} \bigl(\bE S_t(B)
\bigr)^2 .
\end{eqnarray*}
\end{theorem}

\begin{pf}
We define the point processes
\[
\xi_t:=\frac{1}{k!}\sum_{(x_1,\ldots,x_k)\in\mu_{t,\neq
}^k}
\delta_{f_t(x_1,\ldots,x_k)} , \qquad t\geq1 , %
\]
and denote their intensity measures by $\bL_t$, $t\geq1$. By
construction, $S_t(B)$ and $\xi_t(B)$ follow the same distribution. We
notice that for any fixed $h\in\operatorname{Lip}(1)$ (recall that
these are
all $h\dvtx\RR\to\RR$ whose Lipschitz constant is at most one) and
$B\in
\cY$ the
mapping $\o\mapsto h(\o(B))$ from $\widetilde{\N}_\YY$ to $\RR$ satisfies
\[
\bigl\llvert h\bigl(\o_1(B)\bigr)-h\bigl(\o_2(B)\bigr)
\bigr\rrvert\leq\bigl\llvert\o_1(B)-\o_2(B)\bigr\rrvert
\leq\dTV(\o_1,\o_2), \qquad\o_1,
\o_2\in\widetilde{\N}_\YY, %
\]
and thus belongs to $\cC$. Consequently, if $\zeta_t$ is a
Poisson process on $\YY$ with intensity measure $\bL_t$, the
definition of the Wasserstein distance and \eqref
{eqRubinsteinDistanceEquivalent} yield
\begin{eqnarray*}
\dW\bigl(S_t(B),\zeta_t(B)\bigr) &=& \dW\bigl(
\xi_t(B),\zeta_t(B)\bigr) = \sup_{h\in
\operatorname{Lip}(1)}
\bigl\llvert\bE h\bigl(\xi_t(B)\bigr)-\bE h\bigl(
\zeta_t(B)\bigr)\bigr\rrvert
\\
&\leq&\sup_{g\in\cC}\bigl\llvert\bE g(\xi_t\mid
_B)-\bE g(\zeta_t\mid_B)\bigr\rrvert=
\dKR(\xi_t\mid_B,\zeta_t\mid
_B) .
\end{eqnarray*}
Now Theorem~\ref{thmMain} and the observation that $\bL_t(B)=\bE
S_t(B)$ imply the result for the choice $\lambda=\bE S_t(B)$. The
general case follows from the triangle inequality for the Wasserstein
distance and the fact that the Wasserstein distance between a Poisson
random variable with mean $\bE S_t(B)$ and another Poisson random
variable with mean $\lambda$ is bounded by $\llvert\bE S_t(B)-\lambda
\rrvert$.
\end{pf}

We emphasize that Theorem~\ref{thmUstatistics} deals with Poisson
approximation in Wasserstein distance. As already stated in \eqref
{eqInequalityDistances}, this is stronger than approximation in total
variation distance, which is usually considered in the literature; see
\cite{BarbourXia2006} for the only exception we are aware of. This is
possible thanks to our functional limit Theorem~\ref{thmMain}, which
deals with the Kantorovich--Rubinstein distance rather than the total
variation distance for point processes.

The Poisson approximation in total variation distance of U-statistics
over binomial input was considered in \cite{BarbourEagleson}. If we
assume that $\bE S_t(B)=\lambda$ for $t\geq1$ for the binomial case
in Theorem~\ref{thmUstatistics}, we obtain up to a constant, which
may depend on~$\lambda$, the same bound as in \cite{BarbourEagleson},
Theorem~2, for the total variation distance.

In \cite{Peccati11ChenStein}, an abstract bound for the Poisson
approximation of Poisson functionals (i.e., random variables depending
on a Poisson process) is derived, which is also applicable to
U-statistics over Poisson input. Our Theorem~\ref{thmUstatistics}
yields better rates of convergence for this special class of Poisson
functionals. In fact, the bound in~\cite{STScalingLimits},
Proposition~4.1, which is derived from \cite
{Peccati11ChenStein}, involves the square root of $\hat{r}_t(B)$ [see
Remark~\ref{remnonUniform}(iii)], while in the bound for the Poisson
case in Theorem~\ref{thmUstatistics} only $\hat{r}_t(B)$ enters.


To illustrate the use of Theorem~\ref{thmUstatistics} let us consider
a particular example, which will recur also in the following
subsections. Let $K\subset\RR^d$ ($d\geq1$) be a compact convex set
with volume
one. For $t\geq1$ let $\eta_t$ be a homogeneous Poisson process in
$K$ of intensity $t$, and denote by $\beta_t$ a binomial process in
$K$ with $\lceil
t\rceil$ points distributed according to the uniform distribution on
$K$. For a family $(\theta_t)_{t\geq1}$ of positive real numbers let
us construct
the random geometric graph with vertex set~$\mu_t$, where $\mu_t$ is
$\eta_t$ or $\beta_t$, by drawing an edge between two distinct
vertices $y_1$ and $y_2$ whenever their Euclidean distance $\llVert
y_1-y_2\rrVert
$ is bounded by $\theta_t$. These random graphs are the natural
geometric counterparts to the
classical Erd\"os--R\'enyi models for combinatorial random graphs. For
background material we refer the reader to the monograph
\cite{PenroseBook} and also to the recent paper \cite{RSTGilbert} as well
as the references cited therein.

For the random geometric graph introduced above, let $E_t$ be the
number of edges. Note that $E_t$ is a U-statistic of the
form
\[
E_t=\frac{1}{2}\sum_{(y_1,y_2)\in\mu_{t,\neq}^2}{\mathbf1}
\bigl(\llVert y_1-y_2\rrVert\leq\theta_t
\bigr) . %
\]
The multivariate Mecke formula \eqref{eqMecke} and a computation
using spherical coordinates show that $E_t$ has expectation $t^2
(\k_d\theta_t^d+O(\theta_t^{d+1}) )/2$ in the Poisson case, as
$\theta_t\to0$. For an\vspace*{1pt} underlying binomial process the expected
number of edges is $\lceil t\rceil(\lceil t\rceil-1) (\k_d\theta
_t^d+O(\theta_t^{d+1}) )/2$, as $\theta_t\to0$. If the
expectation of $E_t$ converges to a constant, as $t\to\infty$, $E_t$
can be well approximated by a Poisson random variable. In contrast to
\cite{Peccati11ChenStein}, Theorem 5.1, whose proof involves various
nontrivial computations, we can deduce a corresponding approximation
result from Theorem~\ref{thmUstatistics}; the proof is postponed to
Section~\ref{secGilbertGraph}.

%
\begin{corollary}\label{corPoissonGilbert}
Assume that $\lim_{t\to\infty}t^2\theta_t^d=\lambda\in[0,\infty
)$, and let $Z$ be a Poisson distributed random variable with mean
$\kappa_d\lambda/2$. Then there is a constant $c>0$ only depending on
the space dimension $d$, the set $K$ and $\sup_{t\geq1}t^2\theta
_t^d$ such that
\[
\dW(E_t,Z)\leq c \bigl(\bigl\llvert t^2
\theta_t^d-\lambda\bigr\rrvert+t^{-\min\{2/d,1\}
} \bigr),
\qquad t\geq1 . %
\]
\end{corollary}

%
\begin{remark}
Using the classical Chen--Stein method for Poisson approximation,
Theorem 3.4 in \cite{PenroseBook} delivers a version of Corollary~\ref
{corPoissonGilbert} with the same rate of convergence in the total
variation distance in case of an underlying binomial process. For the
Poisson case, Theorem 3.12(iii) in \cite{PeccatiLachiezeRey} is a
qualitative version of Corollary~\ref{corPoissonGilbert}, which has
been established by the method of moments, and Theorem 5.1 in \cite
{Peccati11ChenStein} adds a total variation bound.
Corollary~\ref{corPoissonGilbert} extends these results to a stronger
probability metric and at the same time improves the rates of
convergence in \cite{Peccati11ChenStein}. Namely, for space dimensions
$d\in\{1,2\}$, Corollary~\ref{corPoissonGilbert} yields an upper
bound of order $\llvert t^2\theta_t^d-\lambda\rrvert+t^{-1}$ (for
the Wasserstein
distance), while Theorem 5.1 in \cite{Peccati11ChenStein} delivers an
upper bound of order $\llvert t^2\theta_t^d-\lambda\rrvert+t^{-1/2}$
(for the
total variation distance).
\end{remark}

\subsection{Compound Poisson approximation of U-statistics}\label{secCompoundPoisson}

As in the previous subsection, we denote by $\mu_t$, $t\geq1$, a
Poisson process $\eta_t$ or a binomial process $\beta_t$ on a lcscH
space $\XX$. For $k\in\NN$ and measurable functions $h_t\dvtx\XX^k\to
\RR$, $t\geq1$, we consider the family of U-statistics
\[
S_t:=\frac{1}{k!}\sum_{(x_1,\ldots,x_k)\in\mu^k_{t,\neq}}
h_t(x_1,\ldots,x_k) , \qquad t\geq1 .
\]
Since the sum runs also over all permutations of a fixed $(x_1,\ldots
,x_k)\in\mu_{t,\neq}^k$, we assume without loss of generality that
$h_t$ is symmetric for any $t\geq1$. For a fixed constant $\gamma\in
\RR$ and $t\geq1$, we define
\[
\bL_t(A):=\frac{1}{k!}\bE\sum_{(x_1,\ldots,x_k)\in\mu^k_{t,\neq
}}
{\mathbf1}\bigl(h_t(x_1,\ldots,x_k)\in
t^{-\gamma} A\setminus\{0\}\bigr) , \qquad A\in\cB(\RR) , %
\]
and
\begin{eqnarray*}
r_t&:=&\max_{1\leq\ell\leq k-1} \int_{\XX^\ell}
\biggl( \int_{\XX
^{k-\ell}} {\mathbf1}\bigl(h_t(x_1,
\ldots,x_k)\neq0\bigr) \bK_t^{k-\ell
}\bigl(\dint(x_{\ell+1},\ldots,x_k)\bigr) \biggr)^2
\\
&&\hspace*{50pt}{}\times \bK_t^\ell\bigl(\dint(x_1,\ldots,x_\ell)
\bigr) %
\end{eqnarray*}
for $k\geq2$, and put $r_t:=0$ if $k=1$. The following result compares
the U-statistic $S_t$ with a compound Poisson random variable. Most of
the existing literature is based on a direct use of Stein's method, but
only for discrete compound Poisson random variables. This approach is
technically sophisticated and also needs, in general, certain
monotonicity assumptions. Moreover, there are even situations in which
the solution of the so-called Stein equation cannot be controlled
appropriately, and hence in which Stein's method is of little use; see
\cite{BarbourUtev1}. Being a consequence of the functional limit
theorem (Theorem~\ref{thmMain}), our approach circumvents such
technicalities and also allows us to deal with compound Poisson random
variables having a discrete or continuous distribution.

%
\begin{theorem}\label{thmCompoundPoisson}
Let $\zeta$ be a Poisson process on $\RR$ with a finite intensity
measure $\bM$, let $Z:=\sum_{x\in\zeta} x$
and let $\gamma\in\RR$. Then
\[
\dTVRandom\bigl(t^\gamma S_t,Z\bigr)\leq\dTV(
\bL_t,\bM)+\frac{2^{k+1}}{k!} r_t , \qquad t\geq1
\]
if in the definition of $S_t$ a Poisson process $\eta_t$ is used, and
\[
\dTVRandom\bigl(t^\gamma S_t,Z\bigr)\leq\dTV(
\bL_t,\bM)+\frac{2^{k+1}}{k!} r_t+\frac{6^k k!}{t}
\bL_t(\RR)^2 , \qquad t\geq1 %
\]
if the underlying point process is a binomial process $\beta_t$.
\end{theorem}

\begin{pf}
We consider the point processes
\[
t^\gamma\bullet\xi_t:=\frac{1}{k!}\sum
_{(x_1,\ldots,x_k)\in\mu
^k_{t,\neq}} {\mathbf1}\bigl(h_t(x_1,
\ldots,x_k)\neq0\bigr) \delta_{t^\gamma
h_t(x_1,\ldots,x_k)} , \qquad t\geq1 .
\]
It follows from the definition of the total variation distance and
\eqref{eqRubinsteinDistanceEquivalent} that
\[
\dTVRandom\bigl(t^\gamma S_t, Z\bigr)=\sup
_{A\in\mathcal{B}(\RR)}\biggl\llvert\bE{\mathbf1} \biggl(\sum
_{x\in t^\gamma\bullet\xi_t}x\in A \biggr) - \bE{\mathbf1} \biggl(\sum
_{x\in\zeta} x \in A \biggr)\biggr\rrvert\leq\dKR
\bigl(t^\gamma\bullet\xi_t, \zeta\bigr) %
\]
since the maps $\o\to{\mathbf1}(\sum_{x\in\o}x\in A)$ belong to $\cC
$. Now Theorem~\ref{thmMain} implies that
\[
\dKR\bigl(t^\gamma\bullet\xi_t,\zeta\bigr) \leq\dTV(
\bL_t,\bM)+\frac
{2^{k+1}}{k!} r_t , \qquad t\geq1 ,
\]
and
\[
\dKR\bigl(t^\gamma\bullet\xi_t,\zeta\bigr) \leq\dTV(
\bL_t,\bM)+\frac
{2^{k+1}}{k!} r_t+\frac{6^k k!}{t}
\bL_t(\RR)^2 , \qquad t\geq1 , %
\]
for the Poisson and the binomial case, respectively. This completes the proof.
\end{pf}

%
\begin{remark}
A compound Poisson random variable $Z$ can alternatively be written as
$Z=\sum_{i=1}^N X_i$, where $N$ is a Poisson distributed random
variable and $(X_i)_{i\in\NN}$ is a sequence of independent and
identically distributed random variables such that $N$ and $(X_i)_{i\in
\NN}$ are independent. However, the representation of $Z$ in terms of
the Poisson process $\zeta$ fits better into our general framework.
\end{remark}

For the compound-Poisson approximation of U-statistics in the binomial
case, a bound similar to that in Theorem~\ref{thmCompoundPoisson} is
derived in \cite{EichelsbacherRoos}, Section~3.6. However, in that
paper $h_t$ is required to take values in the nonnegative integers,
whereas we do not need to impose such a condition. In addition, we are
not aware of any analogous result for an underlying Poisson process.

As an application of Theorem~\ref{thmCompoundPoisson} we consider
general edge-length functionals of the random geometric graph
introduced in the course of the previous subsection. Fix a parameter
$b\in\RR$, and define
\[
L^{(b)}_t:=\frac{1}{2} \sum
_{(x_1,x_2)\in\mu^2_{t,\neq}} {\mathbf1}\bigl(\dist(x_1,x_2)\leq
\theta_t\bigr) \dist(x_1,x_2)^b ,
\qquad t\geq1 , %
\]
where $\mu_t$ stands either for a Poisson process $\eta_t$ or a
binomial process $\beta_t$.
In particular, $L^{(0)}_t$ is the number of edges in the random
geometric graph, and $L_t^{(1)}$ is its total edge length. As in
Section~\ref{secUstatistics}, we consider the situation where the
distance parameters $(\theta_t)_{t\geq1}$ are chosen in such a way
that the expected number of edges converges to a constant, as $t\to
\infty$. Recall that in Corollary~\ref{corPoissonGilbert} the number
of edges $L_t^{(0)}$ has been approximated by a Poisson random
variable. For general exponents $b$ we approximate $L_t^{(b)}$ by a
suitable compound Poisson random variable. The proof of the next result
is postponed to Section~\ref{secGilbertGraph} below.

%
\begin{corollary}\label{corCompoundPoissonGilbert}
Fix $b\in\RR$, and assume that $\lim_{t\to\infty}t^2\theta
_t^d=\lambda\in[0,\infty)$. Define $Z:=\sum_{i=1}^N \llVert
X_i\rrVert^b$,
where $N$ is a Poisson distributed random variable with mean $\kappa
_d\lambda/2$ and $(X_i)_{i\in\NN}$ are independent and uniformly
distributed points in $B^d(\lambda^{1/d})$, which are independent of
$N$. Then there is a constant $c>0$ only depending on the space
dimension $d$, the set $K$ and $\sup_{t\geq1}t^2\theta_t^d$ such that
\[
\dTVRandom\bigl(t^{2b/d}L_t^{(b)},Z\bigr)\leq c
\bigl(\bigl\llvert t^2\theta_t^d-\lambda\bigr
\rrvert+t^{-\min\{2/d,1\}} \bigr), \qquad t\geq1 . %
\]
\end{corollary}

%
\begin{remark}
Corollary~\ref{corCompoundPoissonGilbert} without a rate of
convergence has been derived in \cite{RSTGilbert}, Theorem 3.5, by
combining a point process convergence result with the continuous
mapping theorem. Thanks to Theorem~\ref{thmCompoundPoisson} we are
able to add a rate of convergence for the total variation distance.
\end{remark}

\subsection{Approximation of U-statistics by \texorpdfstring{$\alpha$}{alpha}-stable random variables}\label{secstable}

Let us denote by $\mu_t$, $t\geq1$, a Poisson process $\eta_t$ or a
binomial process $\beta_t$ as in the previous subsections. For fixed
$k\in\NN$ and measurable functions $h_t\dvtx \XX^k\to\RR$, $t\geq1$, let
\[
S_t:= \frac{1}{k!}\sum_{(x_1,\ldots,x_k)\in\mu^k_{t,\neq}}
h_t(x_1,\ldots,x_k) , \qquad t\geq1 .
\]
Here, we can and will assume without loss of generality that $h_t$ is
symmetric for any $t\geq1$. We are interested in the limiting behavior
of these U-statistics in situations where their summands are heavy
tailed, and approximate $S_t$ by an $\alpha$-stable random variable
$Z$. Recall that this means that for any $n\in\NN$ there are
independent copies $Z_1,\ldots,Z_n$ of $Z$ satisfying the
distributional equality
$n^{-1/\a}(Z_1+\cdots+Z_n)\stackrel{D}{=}Z$. We fix $\alpha\in
(0,1)$ and $\gamma\in\RR$ and apply our functional limit theorem to
the point processes
%
\begin{eqnarray}
t^\gamma\bullet\xi_t:=\frac{1}{k!}\sum
_{(x_1,\ldots,x_k)\in\mu
^k_{t,\neq}} {\mathbf1}\bigl(h_t(x_1,
\ldots,x_k)\neq0\bigr) \delta_{\operatorname{sign}(h_t(x_1,\ldots
,x_k)) t^\gamma\llvert h_t(x_1,\ldots
,x_k)\rrvert^{-\alpha}} ,\nonumber
\\[-4pt]
\eqntext{t\geq1 ,}
\end{eqnarray}
on $\RR$, where $\operatorname{sign}(a)={\mathbf1}(a\geq0)-{\mathbf
1}(a<0)$. If $\mu_t$ is a binomial process, the convergence of the
U-statistic $S_t$ to an $\alpha$-stable random variable was considered
in~\cite{DehlingEtAl} without giving rates of convergence. Thanks to
our quantitative bound for the Kantorovich--Rubinstein distance in
Theorem~\ref{thmMain}, we are in the position to add a rate of
convergence for the Kolmogorov distance. The statement of our result is
prepared by introducing some notation. For $A\in\cB(\RR)$ and $t\geq
1$, we define
\begin{eqnarray*}
\bL_t(A) &:=& \frac{1}{k!}\bE\sum_{(x_1,\ldots,x_k)\in\mu^k_{t,\neq}}
{\mathbf1}\bigl(h_t(x_1,\ldots,x_k)\neq0\bigr)
\\
&&\hspace*{75pt}{}\times {\mathbf1}\bigl(\operatorname{sign}\bigl(h_t(x_1,
\ldots,x_k)\bigr) \bigl\llvert h_t(x_1,
\ldots,x_k)\bigr\rrvert^{-\alpha}\in t^{-\gamma}A\bigr),
\end{eqnarray*}
which is the intensity measure of $t^\gamma\bullet\xi_t$, and
\begin{eqnarray*}
r_t(A) & :=& \max_{1\leq\ell\leq k-1}\int_{\XX^\ell}
\biggl( \int_{\XX^{k-\ell}} {\mathbf1}\bigl(h_t(x_1,
\ldots,x_k)\neq0\bigr)
\\
&&\hspace*{82pt}{}\times  {\mathbf1}\bigl(\operatorname{sign}
\bigl(h_t(x_1,\ldots,x_k)\bigr) \bigl\llvert
h_t(x_1,\ldots,x_k)\bigr\rrvert
^{-\alpha}\in t^{-\gamma}A\bigr)
\\
&&\hspace*{195pt}{}\times  \bK_t^{k-\ell}\bigl(\dint(x_{\ell+1},
\ldots,x_k)\bigr) \biggr)^2
\\
&&\hspace*{51pt}{}\times \bK_t^\ell \bigl(\dint(x_1,\ldots,x_\ell)\bigr)
\end{eqnarray*}
if $k\geq2$ and $r_t(A):=0$ if $k=1$. The following result contains a
quantitative bound for the approximation of U-statistics by an $\alpha
$-stable random variable with $\alpha\in(0,1)$. 

%
\begin{theorem}\label{thmstable}
Let $\alpha\in(0,1)$, and let $\bM$ be either the Lebesgue measure
on~$\RR$ or its restriction to $\RR_+$. Define $Z:=\sum_{x\in\zeta
} \operatorname{sign}(x) \llvert x\rrvert^{-1/\alpha}$, where $\zeta
$ is a
Poisson process with intensity measure $\bM$. Assume that there are a
constant $\gamma\in\RR$ and functions $g_1,g_2,g_3\dvtx \RR_+^2\to\RR
_+$ such that, for any $a>0$ and $t\geq1$,
%
\begin{equation}
\label{eqbounddTV} \dTV\bigl(\bL_t\llvert_{[-a,a]},\bM\rrvert
_{[-a,a]}\bigr) \leq g_1(a,t) , \qquad r_t
\bigl([-a,a]\bigr)\leq g_2(a,t)
\end{equation}
and
%
\begin{eqnarray}\label{eqboundRestSt}
&& \frac{t^{-\gamma/\alpha}}{k!} \bE\sum
_{(x_1,\ldots,x_k)\in\mu
^k_{t,\neq}} {
\mathbf1}\bigl(\bigl\llvert h_t(x_1,\ldots,x_k)
\bigr\rrvert< t^{\gamma/\alpha
}a^{-1/\alpha} \bigr) \bigl\llvert
h_t(x_1,\ldots,x_k)\bigr\rrvert
\nonumber\\[-8pt]\\[-8pt]\nonumber
&&\qquad \leq g_3(a,t) .
\end{eqnarray}
Then there is a constant $C>0$ only depending on $\alpha$ and $k$ such that
\[
\dK\bigl(t^{-\gamma/\alpha}S_t,Z\bigr)\leq C g(t) , \qquad t\geq1 ,
\]
where
\begin{eqnarray*}
g(t):=\cases{\displaystyle\inf_{a>0}\max\bigl\{a^{1/2-1/(2\alpha)},g_1(a,t),g_2(a,t),
\sqrt{g_3(a,t)}\bigr\}, &\quad$\mu_t= \eta_t$,
\vspace*{3pt}\cr
\displaystyle \inf_{a>0}\max\bigl \{a^{1/2-1/(2\alpha)},g_1(a,t),g_2(a,t),\sqrt{g_3(a,t)},a^2/t\bigr\}, &\quad$\mu_t=
\beta_t$.}
\end{eqnarray*}
\end{theorem}

\begin{pf}
For $a>0$ we define the random variables
\[
S_{t,a}:=\frac{1}{k!}\sum_{(x_1,\ldots,x_k)\in\mu^k_{t,\neq}} {
\mathbf1}\bigl(\bigl\llvert h_t(x_1,\ldots,x_k)
\bigr\rrvert\geq t^{\gamma/\alpha}a^{-1/\alpha}\bigr) h_t(x_1,
\ldots,x_k) , \qquad t\geq1 , %
\]
and
\[
Z_a:=\sum_{x\in\zeta} {\mathbf1}\bigl(\llvert x
\rrvert\leq a\bigr) \operatorname{sign}(x) \llvert x\rrvert^{-1/\alpha
} .
\]
Then, for any $a>0$ and $\varepsilon>0$, we find that
\begin{eqnarray*}
&& \dK\bigl(t^{-\gamma/\alpha}S_t,Z\bigr)
\\
&&\qquad  \leq \bP\bigl(t^{-\gamma/\alpha}
\llvert S_t-S_{t,a}\rrvert\geq\varepsilon\bigr) +\dK
\bigl(t^{-\gamma/\alpha
}S_{t,a},Z\bigr)
\\
&&\quad\qquad{} +\sup_{z\in\RR}\bigl
\llvert\bP(Z\leq z)-\bP(Z\leq z+\varepsilon)\bigr\rrvert
\\
&&\qquad \leq\bP\bigl(t^{-\gamma/\alpha}\llvert S_t-S_{t,a}\rrvert
\geq\varepsilon\bigr)+\bP\bigl(\llvert Z-Z_a\rrvert\geq\varepsilon
\bigr) +\dK\bigl(t^{-\gamma/\alpha}S_{t,a},Z_a\bigr)
\\
&&\quad\qquad{} +2\sup_{z\in\RR}\bigl\llvert\bP(Z\leq z)-\bP(Z\leq z+
\varepsilon)\bigr\rrvert.
\end{eqnarray*}
Combining Markov's inequality with the multivariate Mecke formula
\eqref{eqMecke} and assumption \eqref{eqboundRestSt}, we obtain
that, for all $\varepsilon>0$,
\[
\bP\bigl(\llvert Z-Z_a\rrvert\geq\varepsilon\bigr) \leq
\frac{2}{\varepsilon} \int_a^\infty x^{-1/\alpha}\,\dint x =\frac{2 a^{1-1/\alpha
}}{(1/\alpha-1)\varepsilon}
\]
and
\[
\bP\bigl(t^{-\gamma
/\alpha}\llvert
S_t-S_{t,a}\rrvert\geq\varepsilon\bigr) \leq
\frac
{g_3(a,t)}{\varepsilon} . %
\]
As $\alpha$-stable random variable, $Z$ has a bounded density; see
\cite{Zolotarev}, page 13. Hence there is a constant $C_\alpha>0$
only depending on $\alpha$ such that
\[
\sup_{z\in\RR}\bigl\llvert\bP(Z\leq z)-\bP(Z\leq z+\varepsilon)
\bigr\rrvert\leq C_\alpha\varepsilon, \qquad\varepsilon\geq0 .
\]
It follows from the definition of the Kolmogorov distance and \eqref
{eqRubinsteinDistanceEquivalent} that
\begin{eqnarray*}
\dK\bigl(t^{-\gamma/\alpha}S_{t,a},Z_a\bigr) &=&\sup
_{z\in\RR} \biggl\llvert\bP\biggl(\sum
_{x\in t^\gamma\bullet\xi_t} {\mathbf1}\bigl(x\in[-a,a]\bigr)
\operatorname{sign}(x)
\llvert x\rrvert^{-1/\alpha}\leq z \biggr)
\\
&&\hspace*{20pt}{}  - \bP\biggl(\sum
_{x\in\zeta} {\mathbf1}\bigl(x\in[-a,a]\bigr) \operatorname{sign}(x)
\llvert x\rrvert^{-1/\alpha} \leq z \biggr)\biggr\rrvert
\\
& \leq&\dKR\bigl(t^{\gamma}\bullet\xi_t\mid_{[-a,a]},
\zeta\mid_{[-a,a]}\bigr) .
\end{eqnarray*}
Now we consider the Poisson case and the binomial case separately. For
an underlying Poisson process, Theorem~\ref{thmMain} and the
assumptions in \eqref{eqbounddTV} show that
\[
\dKR\bigl(t^{\gamma}\bullet\xi_t\mid_{[-a,a]},\zeta
\mid_{[-a,a]}\bigr) \leq g_1(a,t) + \frac{2^{k+1}}{k!}
g_2(a,t) , \qquad t\geq1 . %
\]
Combining this with the previous estimates, we see that
\[
\dK\bigl(t^{-\gamma/\alpha}S_t,Z\bigr) \leq\frac{2 a^{1-1/\alpha
}}{(1/\alpha-1)\varepsilon}+
\frac{g_3(a,t)}{\varepsilon}+2C_\alpha\varepsilon+g_1(a,t) +
\frac{2^{k+1}}{k!} g_2(a,t) . %
\]
Thus choosing $\varepsilon= \sqrt{\max\{a^{1-1/\alpha},g_3(a,t)\}}$
yields the assertion. For the binomial case, Theorem~\ref{thmMain}
and the assumptions in \eqref{eqbounddTV} imply that
\begin{eqnarray*}
\dK(S_{t,a},Z_a) &\leq& \frac{2 a^{1-1/\alpha}}{(1/\alpha
-1)\varepsilon}+
\frac{g_3(a,t)}{\varepsilon}+2C_\alpha\varepsilon+g_1(a,t)
\\
&&{} + \frac{2^{k+1}}{k!} g_2(a,t)+\frac{6^k
k!}{t}\bigl(8a^2+2g_1(a,t)^2
\bigr) ,
\end{eqnarray*}
where we have used that $\bL_t([-a,a])^2 \leq(2a+g_1(a,t))^2\leq
8a^2+2g_1(a,t)^2$. Now the same choice for $\varepsilon$ as in the
Poisson case and the fact that the Kolmogorov distance is bounded by
one complete the proof.
\end{pf}

%
\begin{remark}\label{remAlpha>=1}
For all choices of $\alpha\in(0,2]$ there are $\alpha$-stable random
variables, and one can think of U-statistics converging to such
variables. For $\alpha\in(1,2]$ and the binomial case this problem
was considered in \cite{DehlingEtAl,HeinrichWolf1993,MalevichAbad}. A
technique similar to that used in the proof of Theorem~\ref
{thmstable} should also be applicable if $\alpha\in(1,2]$. In this
case the limiting random variable is given by $Z:=\lim_{a\to\infty}
Z_a-\bE Z_a$, whence an additional centering is necessary. In order to
derive bounds similar to those of Theorem~\ref{thmstable}, one has to
control the distance between $Z$ and $Z_a$, which might be difficult to
tackle. We would like to mention that the bounds derived in \cite
{HeinrichWolf1993} also involve a quantity similar to $\dK(Z,Z_a)$.
\end{remark}

To give an application of Theorem~\ref{thmstable}, let us consider
the following distance-power statistics, which are closely related to
the edge functionals of random geometric graphs considered above. Let
for some $d\geq1$, $K\subset\RR^d$ be a compact convex set with
volume one, and let $\bK$ be the restriction of the Lebesgue measure
to $K$. Let $\eta_t$ be a Poisson process in $K$ with intensity
measure $\bK_t=t\bK$, $t\geq1$, and let $\beta_t$, $t\geq1$ be a
binomial process of $\lceil t \rceil$ points, which are independent
and uniformly distributed in $K$. Our aim is to investigate the
limiting behavior of the U-statistics
\[
S_t:=\frac{1}{2}\sum_{(x_1,x_2)\in\mu^k_{t,\neq}} \dist
(x_1,x_2)^{-\tau} , \qquad t\geq1 , %
\]
where $\tau>0$ and $\mu_t$ stands for $\eta_t$ or $\beta_t$. The
following result, whose proof will be given in Section~\ref
{secGilbertGraph} below, deals with the case $\tau>d$.

%
\begin{corollary}\label{corStableGilbert}
Let $\tau>d$, let $\zeta$ be a homogeneous Poisson process on $\RR
_+$ with intensity one and let $Z:=(\kappa_d/2)^{\tau/d}\sum_{x\in
\zeta} x^{-\tau/d}$. Then there is a constant $C>0$ only depending on
$K$, $\tau$ and $d$ such that
\[
\dK\bigl(t^{-2\tau/d} S_t,Z\bigr) \leq C t^{\varrho} ,
\qquad t\geq1 , %
\]
with
\[
\varrho:=\inf_{u>0}\max\biggl\{\frac{1}{2}u-
\frac{\tau}{2d}u, 2u-1, u+\frac{1}{d}u-\frac{2}{d}\biggr\} .
\]
\end{corollary}

%
\begin{example}
To have a more specific example, take $\tau=2d$ in Corollary~\ref
{corCompoundPoissonGilbert}, in which case $\varrho$ has the form
\[
\varrho=\inf_{u>0}\max\biggl\{-{u\over2},2u-1,u+
{u-2\over d}\biggr\} . %
\]
For $d\in\{1,2\}$ the infimum is attained at $u={2\over5}$, giving
that $\varrho=-{1\over5}$. For $d\geq3$, the infimum is attained at
$u=\frac{4}{3d+2}$ so that $\varrho=-\frac{2}{3d+2}$ in this case. Thus
\[
\dK\bigl(t^{-4}S_t,Z\bigr)\leq\cases{ C t^{-1/5},
&\quad$d\in\{1,2\}$,
\vspace*{3pt}\cr
C t^{-2/(3d+2)}, &\quad$d\geq3$,} %
\]
where the $1/2$-stable random variable $Z$ is of the form $Z=c_d\sum
_{x\in\zeta}x^{-2}$ for a unit-intensity homogeneous Poisson process
$\zeta$ on $\RR_+$ and with $c_d=\k_d^2/4$. The distribution of $Z$
can be characterized more explicitly. Namely, applying \cite
{KallenbergFoundations}, Lemma~12.2(i), we see that for all $t\in\RR$,
\begin{eqnarray*}
\bE\exp(\mathfrak{i}tZ) &=& \bE\exp\biggl(\mathfrak{i}t c_d\sum
_{x\in\zeta}x^{-2} \biggr) = \exp\biggl(\int
_0^\infty\bigl(e^{\mathfrak
{i}tc_dx^{-2}}-1\bigr)\,\dint x
\biggr)
\\
&=& \exp(-\sqrt{-\mathfrak{i}t\pi c_d} ) ,
\end{eqnarray*}
where $\mathfrak{i}$ is the imaginary unit. This is the characteristic
function of a centred L\'evy distribution with scale parameter $\pi
c_d/2$. Thus $Z$ has density $x\mapsto{1\over2}\sqrt{c_d/x^3} \exp
(-\pi c_d/(4x)) {\mathbf1}(x>0)$.
\end{example}

%
\begin{remark}
Note that if $\tau<d/2$, then $S_t$ satisfies a central limit theorem
as shown in Theorem 3.1 of \cite{RSTGilbert}. Moreover, the choice
$d/2\leq\tau\leq d$ corresponds to the situation $\alpha\in[1,2]$,
to which Remark~\ref{remAlpha>=1} applies.
\end{remark}

\subsection{Random geometric graphs}\label{secGilbertGraph}

Let $K\subset\RR^d$ ($d\geq1$) be a compact convex set with volume
one. For $t\geq1$ let $\mu_t$ either be a homogeneous Poisson process
$\eta_t$ of intensity $t\geq1$ in $K$ or a binomial process $\beta
_t$ of $\lceil t\rceil$ independent and uniformly distributed points
in $K$, and let $(\theta_t)_{t\geq1}$ be a family of positive real
numbers. Based on this data we construct a random geometric graph as
explained in Section~\ref{secUstatistics}. In contrast to Corollaries
\ref{corPoissonGilbert}~and~\ref
{corCompoundPoissonGilbert}, where $\lim_{t\to\infty}t^2\theta
_t^d=\lambda\in[0,\infty)$, we assume at first that $\lim_{t\to
\infty}t^2\theta_t^d=\infty$ and are interested in the point process
$\xi_{t,a}$ on $K$ defined by
\[
\xi_{t,a}:=\frac{1}{2} \sum_{(x,y)\in\mu_{t,\neq}^2}{
\mathbf1}\bigl(\llVert x-y\rrVert\leq\min\bigl\{\theta_t,t^{-2/d}a
\bigr\}\bigr) \delta_{(x+y)/2} %
\]
for\vspace*{1pt} some $a>0$. In other words, $\xi_{t,a}$ charges the collection of
all midpoints of edges of the random geometric graph whose length does not
exceed $t^{-2/d}a$.

%
\begin{theorem}\label{thmGraphsMultivariate}
Let $a>0$, let $\zeta$ be a Poisson process on $K$ with intensity measure
$\frac{\kappa_d}{2} a^d \mathrm{vol}\mid_K$ and let $\xi_{t,a}$ be
constructed from a Poisson process $\eta_t$ or a
binomial process $\beta_t$ with $t\geq1$. Also suppose that $\lim_{t\to
\infty}t^2\theta_t^d=\infty$. Then $t_0:=\sup\{t\geq1\dvtx
t^2\theta_t^d<a^d\}\cup\{1\}<\infty$, and there is a constant $C>0$
only depending on $a$, $d$ and $K$ such that
\[
\dKR(\xi_{t,a},\zeta)\leq C t^{-\min\{2/d,1\}} , \qquad t>
t_0 . %
\]
\end{theorem}

The rest of this subsection is devoted to the proofs of Theorem~\ref
{thmGraphsMultivariate} as well as Corollaries~\ref
{corPoissonGilbert}, \ref{corCompoundPoissonGilbert} and~\ref{corStableGilbert}. We prepare with the following
lemma. In order to deal with the Poisson and the binomial case in
parallel, we define $\chi(t)=t^2$ and $\widetilde{\chi}(t)=t^3$ if
$\mu_t=\eta_t$ and $\chi(t)=\lceil t\rceil(\lceil t \rceil-1)$ and
$\widetilde{\chi}(t)=(\lceil t\rceil)^3$ if $\mu_t=\beta_t$.

%
\begin{lemma}\label{lemBoundsGilbert}
There is a constant $C_K>0$ only depending on $d$ and $K$ such that
%
\begin{eqnarray}
\label{eqBoundGilbert1} && \biggl\llvert\frac{1}{2}\bE\sum
_{(x,y)\in\mu^2_{t,\neq}} {\mathbf1}\bigl((x+y)/2\in B, \llVert
x-y\rrVert\in
\tilde{A}\bigr)\nonumber
\\
&&\hspace*{32pt}{} - \frac{\kappa_d}{2} \vol(B) t^2 d \int
_0^\infty{\mathbf1}(r\in\tilde{A}) r^{d-1}\,
\dint r \biggr\rrvert
\\
&&\qquad\leq2 C_K \kappa_d t^2 \bigl(
\tilde{a}^{d+1}+\tilde{a}^{2d}\bigr)+\frac{\kappa_d}{2} t
\tilde{a}^d\nonumber
\end{eqnarray}
for all Borel sets $B\subset K$ and $\tilde{A}\subset[0,\tilde{a}]$
with $\tilde{a}>0$. Moreover,
%
\begin{equation}
\label{eqBoundGilbert2} \widetilde{\chi}(t)\int_K \biggl(\int
_K {\mathbf1}\bigl((x+y)/2\in B, \llVert x-y\rrVert\leq u
\bigr)\,\dint x \biggr)^2 \,\dint y \leq8t^3
\kappa_d^2 u^{2d}
\end{equation}
for all Borel sets $B\subset K$ and $u\geq0$.
\end{lemma}

\begin{pf}
By the multivariate Mecke formula \eqref{eqMecke} for the Poisson
process and its analogue \eqref{eqMomentMeasureBinomial} for the
binomial case, we obtain that
\begin{eqnarray*}
&& \frac{1}{2} \bE\sum_{(x,y)\in\mu^2_{t,\neq}} {\mathbf1}
\bigl((x+y)/2\in B, \llVert x-y\rrVert\in\tilde{A}\bigr)
\\
&&\qquad = \frac{\chi(t)}{2} \int_K \int_K
{\mathbf1}\bigl((x+y)/2\in B, \llVert x-y\rrVert\in\tilde{A}\bigr)\,\dint x \,\dint y
\\
&&\qquad = \frac{\chi(t)}{2} \int_{\RR^d} \int_{\RR^d}
{\mathbf1}\bigl((x+y)/2\in B, \llVert x-y\rrVert\in\tilde{A}\bigr)\,\dint x \,\dint y
\\
&&\qquad\quad{}- \frac{\chi(t)}{2} \int_{(\RR^d)^2\setminus
K^2} {\mathbf1}
\bigl((x+y)/2\in B, \llVert x-y\rrVert\in\tilde{A}\bigr)\,\dint(x,y) .
\end{eqnarray*}
To the first term in the last expression we apply
the change of variables $u=x-y$, $v=(x+y)/2$, which has Jacobian one,
and spherical coordinates
to see that
\begin{eqnarray*}
&& \frac{\chi(t)}{2}\int_{\RR^d}\int_{\RR^d}{
\mathbf1} \bigl((x+y)/2\in B, \llVert x-y\rrVert\in\tilde{A} \bigr)\,\dint x\, \dint y
\\
&&\qquad =  \frac{\chi
(t)}{2}\int_{\RR^d}\int_{\RR^d}{
\mathbf1} \bigl(v\in B, \llVert u\rrVert\in\tilde{A} \bigr)\,\dint u\,
\dint v
\\
&&\qquad = \frac{\chi(t)}{2}\vol(B) d\kappa_d \int_0^\infty{
\mathbf1}(r\in\tilde{A}) r^{d-1} \,\dint r .
\end{eqnarray*}
A straightforward compuatation shows that
\[
\biggl\llvert\bigl(t^2-\chi(t)\bigr)\vol(B)\frac{\kappa_d}{2} d
\int_0^\infty{\mathbf1}(r\in\tilde{A})
r^{d-1} \,\dint r \biggr\rrvert\leq\frac{\kappa
_d}{2} t
\tilde{a}^d . %
\]
For the second term we have, independently of $B$, the upper bound
\begin{eqnarray*}
&& \frac{\chi(t)}{2} \int_{(\RR^d)^2\setminus K^2} {\mathbf1} \bigl
((x+y)/2\in B,
\llVert x-y\rrVert\in\tilde{A} \bigr)\,\dint(x,y)
\\
&&\qquad \leq2t^2 \vol\bigl(\bigl\{x\in\RR^d
\setminus K\dvtx \dist(x,K)\leq\tilde{a}\bigr\} \bigr) \kappa_d
\tilde{a}^d .
\end{eqnarray*}
From Steiner's formula \eqref{eqSteiner} it follows that
there is a constant $C_K>0$ only depending on $d$ and $K$ such
that
\[
\vol\bigl(\bigl\{x\in\RR^d\setminus K\dvtx \dist(x,K)\leq\tilde{a}\bigr
\} \bigr)\leq C_{K} \bigl(\tilde{a}+\tilde{a}^d\bigr) .
\]
Combining these estimates yields the first bound. On the other hand, we have
\begin{eqnarray*}
\widetilde{\chi}(t)\int_K \biggl(\int
_K {\mathbf1}\bigl((x+y)/2\in B, \llVert x-y\rrVert \leq u
\bigr)\,\dint x \biggr)^2 \,\dint y &\leq& 8t^3 \int
_K \bigl(\kappa_d u^d
\bigr)^2 \,\dint y
\\
&=& 8t^3 \kappa_d^2 u^{2d} ,
\end{eqnarray*}
which is the second bound.
\end{pf}

\begin{pf*}{Proof of Theorem~\ref{thmGraphsMultivariate}}
Due\vspace*{1pt} to our assumption that
$\lim_{t\to\infty} t^2\theta_t^d=\infty$, we have that $t_0:=\sup
\{t\geq1\dvtx t^2\theta_t^d<a^d\}\cup\{1\}<\infty$. Note that $\min\{
\theta_t,\break t^{-2/d}a\}=t^{-2/d}a$ for $t> t_0$. We denote by $\bL
_{t,a}$ the intensity measure of $\xi_{t,a}$. For $t>t_0$ the choice
$\tilde{A}=[0,\min\{\theta_t,t^{-2/d}a\}]=[0,t^{-2/d}a]$ in \eqref
{eqBoundGilbert1} leads to
\begin{eqnarray*}
&& \biggl\llvert\bL_{t,a}(B)-\frac{\kappa_d}{2}\vol(B) t^2
\bigl(t^{-2/d} a\bigr)^d \biggr\rrvert
\\
&&\qquad \leq2C_K
\kappa_d t^2 \bigl(t^{-2-2/d}a^{d+1}+t^{-4}a^{2d}
\bigr)+\frac{\kappa
_d}{2}t^{-1}a^d %
\end{eqnarray*}
so that $\dTV(\bL_{t,a},\frac{\kappa_d}{2}a^d \mathrm{vol}\mid_K)\leq C_1
t^{-\min\{2/d,1\}}$ for $t>t_0$ with a constant $C_1>0$ only depending
on $a$, $d$ and $K$. Moreover, there is a constant $C_2>0$ only
depending on $a$, $d$ and $K$ such that $\bL_{t,a}(K)\leq C_2$ for all
$t>t_0$. Inequality \eqref{eqBoundGilbert2} implies that for $t>t_0$,
\begin{eqnarray*}
\widetilde{\chi}(t) \int_K \biggl( \int
_K {\mathbf1}\bigl(\llVert x-y\rrVert\leq\min\bigl\{
\theta_t,t^{-2/d}a\bigr\}\bigr)\,\dint x \biggr)^2\,
\dint y &\leq& 8 t^3 \kappa_d^2
\bigl(t^{-2/d}a\bigr)^{2d}
\\
&=& 8 \kappa_d^2 a^{2d} t^{-1}.
\end{eqnarray*}
Now, application of Theorem~\ref{thmMain} completes the proof.
\end{pf*}

\begin{pf*}{Proof of Corollary~\ref{corPoissonGilbert}}
The choice $B=K$ and $\tilde{A}=[0,\theta_t]$ in \eqref
{eqBoundGilbert1} leads to
\begin{eqnarray*}
\biggl\llvert\bE E_t - \frac{\kappa_d}{2}\lambda\biggr\rrvert& \leq&\biggl\llvert\frac
{\kappa_d}{2}\lambda- \frac{\kappa_d}{2}t^2
\theta_t^d\biggr\rrvert+ \biggl\llvert\bE
E_t - \frac{\kappa_d}{2}t^2\theta_t^d
\biggr\rrvert
\\
& \leq&\frac{\kappa_d}{2}\bigl\llvert\lambda-t^2
\theta_t^d\bigr\rrvert+2C_K
\kappa_d t^2\bigl(\theta_t^{d+1}+
\theta_t^{2d}\bigr)+\frac{\kappa_d}{2}t
\theta_t^d
\\
& \leq&\frac{\kappa_d}{2}\bigl\llvert\lambda-t^2
\theta_t^d\bigr\rrvert+2C_K
\kappa_d \biggl(\frac{(\sup_{t\geq1}t^2\theta
_t^d)^{1+1/d}}{t^{2/d}}+\frac
{(\sup_{t\geq1}t^2\theta_t^d)^2}{t^2} \biggr)
\\
&&{}+
\frac{\kappa
_d}{2}\frac{\sup_{t\geq1} t^2\theta_t^d}{t}
\end{eqnarray*}
for $t\geq1$, which also implies that $\bE E_t$ is bounded by a
constant only depending on~$d$, $K$ and $\sup_{t\geq1}t^2\theta_t^d$
for $t\geq1$. It follows from \eqref{eqBoundGilbert2} that
\[
\widetilde{\chi}(t) \int_K \biggl(\int
_K {\mathbf1}\bigl(\llVert x-y\rrVert\leq\theta
_t\bigr)\,\dint x \biggr)^2 \,\dint y \leq8 t^3
\kappa_d^2 \theta_t^{2d}\leq8
\kappa_d^2 \frac{(\sup_{t\geq1}t^2 \theta
_t^d)^2}{t} . %
\]
Now, the assertion is a consequence of Theorem~\ref{thmUstatistics}.
\end{pf*}

\begin{pf*}{Proof of Corollary~\ref{corCompoundPoissonGilbert}}
We assume that $b\neq0$ in the following since for $b=0$ the assertion
follows from Corollary~\ref{corPoissonGilbert}. For a Borel set
$A\subset[0,\infty)$ we define $A^{1/b}:=\{a^{1/b}\dvtx a\in A\setminus
\{
0\}\}$. Hence we have that
\begin{eqnarray*}
\bL_t(A) & :=&\frac{1}{k!} \bE\sum_{(x,y)\in\mu_{t,\neq}^2}
{\mathbf1}\bigl(\llVert x-y\rrVert\leq\theta_t, \llVert x-y\rrVert
^b\in t^{-2b/d} A\setminus\{0\}\bigr)
\\
& =&\frac{1}{k!} \bE\sum_{(x,y)\in\mu_{t,\neq}^2} {\mathbf1}\bigl(
\llVert x-y\rrVert\in t^{-2/d} A^{1/b}\cap[0,
\theta_t]\bigr) .
\end{eqnarray*}
Moreover, we define
\[
\bM(A):=\frac{\kappa_d}{2} d\int_0^{\lambda^{1/d}} {\mathbf1}
\bigl(r\in A^{1/b}\bigr) r^{d-1} \,\dint r , \qquad A\in
\mathcal{B}(\RR) . %
\]
For a Borel set $A\subset[0,\infty)$, inequality \eqref
{eqBoundGilbert1} with $B=K$ and $\tilde{A}=t^{-2/d}A^{1/b}\cap
[0,\theta_t]$ implies that
\begin{eqnarray*}
&& \bigl\llvert\bL_t(A)-\bM(A)\bigr\rrvert
\\
&&\qquad \leq \biggl\llvert
\frac{\kappa_d}{2}t^2 d\int_0^\infty{
\mathbf1}\bigl(r\in t^{-2/d}A^{1/b}\cap[0,\theta_t]\bigr)
r^{d-1} \,\dint r
\\
&&\hspace*{87pt}{} -\frac{\kappa_d}{2} d\int_0^{\lambda^{1/d}}{\mathbf1}\bigl(r\in A^{1/b}\bigr) r^{d-1} \,\dint r \biggr\rrvert
\\
&&\quad\qquad{} + 2C_K\kappa_d t^2\bigl(
\theta_t^{d+1}+\theta_t^{2d}\bigr)+
\frac
{\kappa_d}{2}t \theta_t^d
\\
&&\qquad \leq \frac{\kappa_d}{2}\bigl\llvert\lambda-t^2
\theta_t^d\bigr\rrvert+ 2C_K
\kappa_d \biggl(\frac{(\sup_{t\geq1}t^2\theta
_t^d)^{1+1/d}}{t^{2/d}}+\frac
{(\sup_{t\geq1}t^2\theta_t^d)^2}{t^2} \biggr)
\\
&&\quad\qquad{} +\frac{\kappa_d}{2} \frac{\sup_{t\geq1}t^2\theta_t^d}{t} .
\end{eqnarray*}
Hence there are constants $C_1,C_2>0$ only depending on $d$, $K$ and
$\sup_{t\geq1}t^2\theta_t^d$ such that $\dTV(\bL_t,\bM)\leq C_1
t^{-\min\{2/d,1\}}$ for $t\geq1$ and $\bL_t(\RR)\leq C_2$ for
$t\geq1$. It follows from \eqref{eqBoundGilbert2} that
\[
\widetilde{\chi}(t)\int_K \biggl(\int
_K {\mathbf1}\bigl(\llVert x-y\rrVert\leq\theta
_t\bigr)\,\dint x \biggr)^2 \,\dint y \leq8 t^3
\kappa_d^2 \theta_t^{2d} \leq8
\kappa_d^2 \frac{(\sup_{t\geq1}t^2\theta_t^d)^2}{t} . %
\]
Now, application of Theorem~\ref{thmCompoundPoisson} completes the proof.
\end{pf*}

\begin{pf*}{Proof of Corollary~\ref{corStableGilbert}}
In the following, we check that the assumptions of Theorem~\ref
{thmstable} are satisfied with $h_t(x,y)=(2/\kappa_d)^{\tau/d}\llVert
x-y\rrVert
^{-\tau}$ with $\alpha=d/\tau$ and $\gamma=2$. For a Borel set
$A\subset[0,\infty)$ and $t\geq1$ we have that
\begin{eqnarray*}
\bL_t(A) &:=& \frac{1}{2} \bE\sum_{(x,y)\in\mu^2_{t,\neq}}
{\mathbf1}\bigl(\kappa_d\llVert x-y\rrVert^d/2\in
t^{-2} A\bigr)
\\
&=& \frac{1}{2} \bE\sum_{(x,y)\in
\mu^2_{t,\neq}}
{\mathbf1}\bigl(\llVert x-y\rrVert\in(2/\kappa_d)^{1/d}t^{-2/d}
A^{1/d}\bigr) %
\end{eqnarray*}
with $A^{1/d}:=\{x^{1/d}\dvtx x\in A\}$. In the following let $\bM$ be the
restriction of the Lebesgue measure to $\RR_+$, and let $a>0$. Since
\begin{eqnarray*}
&& \frac{\kappa_d}{2} t^2 d \int_0^\infty
{\mathbf1}\bigl(r\in(2/\kappa_d)^{1/d}t^{-2/d} \bigl(A
\cap[0,a]\bigr)^{1/d}\bigr) r^{d-1} \,\dint r
\\
&&\qquad = \frac{\kappa_d}{2} t^2 \int_0^\infty
{\mathbf1}\bigl(u\in(2/\kappa_d)t^{-2} \bigl(A\cap[0,a]\bigr)
\bigr)\,\dint u= \int_0^\infty{\mathbf1}\bigl(u\in A\cap
[0,a]\bigr)\,\dint u
\\
&&\qquad = \bM\mid_{[0,a]}(A) ,
\end{eqnarray*}
application of \eqref{eqBoundGilbert1} with $B=K$ and $\tilde
{A}=(2/\kappa_d)^{1/d}t^{-2/d} (A\cap[0,a])^{1/d}$ yields that
\begin{eqnarray*}
\bigl\llvert\bL_t\rrvert_{[0,a]}(A)-\bM\mid
_{[0,a]}(A)\bigr\rrvert& \leq&2 C_K\kappa_d
t^2\bigl(c_a^{d+1}t^{-2-2/d}+c_a^{2d}t^{-4}
\bigr)+\frac{\kappa_d}{2} t^{-1}c_a^d
\end{eqnarray*}
with $c_a=(2a/\kappa_d)^{1/d}$. Consequently, there is a constant
$C_1>0$ only depending on $d$ and $K$ such that
\[
\dTV\bigl(\bL_t\llvert_{[0,a]},\bM\rrvert_{[0,a]}
\bigr) \leq C_1 \bigl(a^{1+1/d}t^{-2/d}+a^2
t^{-2}+at^{-1}\bigr)=:g_1(a,t) , \qquad t\geq1 .
\]
It follows from \eqref{eqBoundGilbert2} that
\begin{eqnarray*}
&& \widetilde{\chi}(t)\int_K \biggl(\int
_K {\mathbf1}\bigl(\kappa_d\llVert x-y\rrVert
^d/2\leq t^{-2}a\bigr)\,\dint x \biggr)^2 \,\dint
y
\\
&&\qquad  = \widetilde{\chi}(t) \int_K \biggl(\int
_K {\mathbf1}\bigl(\llVert x-y\rrVert\leq(2/
\kappa_d)^{1/d} t^{-2/d}a^{1/d}\bigr)\,\dint
x \biggr)^2 \,\dint y
\\
&&\qquad \leq8 t^3 \kappa_d^2 (2/
\kappa_d)^2 t^{-4} a^2= 32
t^{-1} a^2=:g_2(a,t) .
\end{eqnarray*}
Moreover, we have that
\begin{eqnarray*}
&&\frac{t^{-2\tau/d}}{2} \bE\sum_{(x,y)\in\mu^2_{t,\neq}} {\mathbf
1}\bigl((2/
\kappa_d)^{\tau/d}\llVert x-y\rrVert^{-\tau}\leq
t^{2\tau/d}a^{-\tau
/d}\bigr) (2/\kappa_d)^{\tau/d}
\llVert x-y\rrVert^{-\tau}
\\
&&\qquad \leq d \kappa_d (2/\kappa_d)^{\tau/d}
t^{2-2\tau
/d}\int_{(\kappa_d/2)^{1/d}t^{-2/d}a^{1/d}}^\infty r^{-\tau}
r^{d-1} \,\dint r
\\
&&\qquad = \frac{d \kappa_d}{\tau-d} (\kappa_d/2)^{1-2\tau/d}
a^{1-\tau/d}=:g_3(a,t) .
\end{eqnarray*}
Now, Theorem~\ref{thmstable} completes the proof.
\end{pf*}

\subsection{Proximity of Poisson flats}\label{secProximity}

For a space dimension $d\geq2$ and a dimension parameter
$m\geq1$ satisfying $m<d/2$, we investigate the mutual arrangement of
the flats of a Poisson $m$-flat process, that is, a Poisson process on
the space of \mbox{$m$-}dimensional affine subspaces of $\RR^d$, which are
called $m$-flats. In order to define such a Poisson $m$-flat process in
a rigorous way, recall that $\GG_m^d$ and $\AA_m^d$
stand for the space of $m$-dimensional linear and
$m$-dimensional affine subspaces of $\RR^d$, respectively. Let $\QQ$
be a
probability measure on $\GG_m^d$ with the property that two
independent random subspaces $L,M\in\GG_m^d$ with distribution $\QQ$
are almost surely in general position, meaning that the dimension of
the linear hull of $L$ and $M$ is $2m$ with probability one.
Note that this is
satisfied, for example, if $\QQ$ is absolutely continuous with respect
to the unique Haar probability measure on $\GG_m^d$; cf. \cite{SW},
Theorem 4.4.5(c). The measure $\QQ$ induces a
translation-invariant measure $\bK_t$ on $\AA_m^d$ via
%
\begin{equation}
\label{eqDefMuFlats} \int_{\AA_m^d}g(E) \bK_t(\dint E)=t
\int_{\GG_m^d}\int_{E_0^\perp
}g(E_0+x)
\vol_{E_0^\perp}(\dint x) \QQ(\dint E_0) ,
\end{equation}
where $t\geq1$ is an intensity parameter, $g\geq0$ is a measurable
function on $\AA_m^d$ and $\vol_{E_0^\perp}$ denotes the Lebesgue
measure on $E_0^\perp$, the orthogonal complement of $E_0$. We use the
convention $\bK:=\bK_1$ and can re-write $\bK_t$ as $\bK_t=t\bK$.
We now consider a Poisson process $\eta_t$
with intensity measure $\bK_t$. This is what is usually called a Poisson
$m$-flat process in stochastic geometry \cite{SW}, Chapter~4.4. One
particular problem for such $m$-flat processes is to describe the
mutual arrangement of the flats in space. Since $m<d/2$, any two
different flats
$E,F$ of $\eta_t$ do not intersect each other with probability
one. Thus they have a well-defined distance $\dist(E,F)$, and we
denote by
$m(E,F)$ the midpoint of the almost surely uniquely determined line
segment realizing this distance (the perpendicular of $E$ and $F$). We
are interested here in the point process of the midpoints $m(E,F)$ such
that the flats $E,F$ are close together, and $m(E,F)$ is in a compact
convex set $K\subset\RR^d$ of volume $0<\vol(K)<\infty$. To the
best of our knowledge, Theorem~\ref{thmFlatsMultivariate} is the
first result describing its asymptotic behavior, as $t\to\infty$. To
do so, we define for $t\geq1$ and $a>0$, $\xi_{t,a}$ on $K$ by
\[
\xi_{t,a}:=\frac{1}{2}\sum_{(E,F)\in\eta_{t,\neq}^2}
\delta_{m(E,F)} {\mathbf1}\bigl(\dist(E,F)\leq at^{-2/(d-2m)},
m(E,F)\in K
\bigr) . %
\]
The intensity measure
$\bL_{t,a}(B)$ of $\xi_{t,a}$ for a Borel set $B\subset K$ is given
by
\[
\bL_{t,a}(B)=\frac{t^2}{2}\int_{\AA_m^d}\int
_{\AA_m^d}{\mathbf1}\bigl(m(E,F)\in B, \dist(E,F)\leq
at^{-2/(d-2m)}\bigr) \bK(\dint E) \bK(\dint F) %
\]
due to the multivariate Mecke formula \eqref{eqMecke}. It follows
from \cite{STFlats}, Theorem 1 (it is readily checked that the
identity there extends from compact convex sets to general Borel sets) that
\[
\bL_{t,a}(B)=\frac{t^2}{2}\kappa_{d-2m}
\bigl(at^{-2/(d-2m)}\bigr)^{d-2m} \vol(B) \int_{\mathbb{G}_k^d}
\int_{\mathbb{G}_k^d} [M,L] \QQ(\dint L) \QQ(\dint M) , %
\]
where $[M,L]$ stands for the subspace determinant of $M$ and $L$
introduced in Section~\ref{secPreliminaries}. This leads to
\[
\bL_{t,a}(B)=\frac{\kappa_{d-2m}}{2} \vol(B) a^{d-2m} \int
_{\GG
_m^d}\int_{\GG_m^d}[L,M] \QQ(\dint L) \QQ(
\dint M) . %
\]
Now, putting
%
\begin{equation}
\label{eqDefBetaFlats} \sC:=\frac{\kappa_{d-2m}}{2} \int_{\GG_m^d}\int
_{\GG_m^d}[L,M] \QQ(\dint L) \QQ(\dint M) ,
\end{equation}
we see that
\[
\dTV\bigl(\bL_{t,a},\sC a^{d-2m} \mathrm{vol}\mid_K
\bigr)=0 , %
\]
where $\mathrm{vol}\mid_K$ stands for the restriction of the Lebesgue measure on
$\RR^d$ to
$K$. Moreover, the proof of \cite{STFlats}, Theorem 3, shows that there
is a constant $\hat{C}>0$ only depending on $a$, $d$, $m$, $\QQ$ and $K$
such that
\[
\hat{r}_t := \sup_{E\in\AA_m^d} t \int
_{\AA_m^d} {\mathbf1}\bigl(m(E,F)\in K, \dist(E,F)\leq
at^{-2/(d-2m)} \bigr) \bK(\dint F)\leq\hat{C} t^{-1} . %
\]
From this we conclude that
\begin{eqnarray*}
r_t &:=& t^3\int_{\AA_m^d} \biggl(\int
_{\AA_m^d} {\mathbf1}\bigl(m(E,F)\in K, \dist(E,F)\leq
at^{-2/(d-2m)} \bigr) \bK(\dint E) \biggr)^2 \bK(\dint F)
\\
&\leq& 2 \hat{C} \bL_{t,a}(K) t^{-1}
\end{eqnarray*}
and in view of Theorem~\ref{thmMain}
the following result for the midpoint process $\xi_{t,a}$.

%
\begin{theorem}\label{thmFlatsMultivariate}
Let $a>0$, and let $\zeta$ be a Poisson process with intensity measure
$\sC a^{d-2m} \mathrm{vol}\mid_K$, where $\sC$ is as at
\eqref{eqDefBetaFlats}. Then there is a constant $C>0$ depending on
$a$, $d$, $m$, $\QQ$ and $K$ such that
\[
\dKR(\xi_{t,a},\zeta)\leq C t^{-1}, \qquad t\geq1 .
\]
\end{theorem}

%
\begin{remark}
(i) Note that because of \eqref{eqIntegralL,M}, the
constant $\sC$ takes the particularly appealing form
\[
\sC=\frac{1}{2}\frac{{d-m\choose m}}{{d\choose m}}\frac{\k
_{d-m}^2}{\k_d} %
\]
if $\QQ$ is the invariant Haar probability measure on $\GG_m^d$ (or,
equivalently, if the $m$-flat process is stationary and isotropic; see
\cite{SW}).

(ii) As opposed to our previous applications, we do not
consider a
binomial counterpart to Theorem~\ref{thmFlatsMultivariate}. The
reason for that is that there is no normalization, which would
turn the measure $\bK_1$ defined at \eqref{eqDefMuFlats} into a
probability measure.

(iii) Theorem~\ref{thmFlatsMultivariate} extends Theorem
\ref{thmGraphsMultivariate} from $m=0$ (which has been excluded
here for technical reasons) to arbitrary $m$ satisfying
$m<d/2$. However, due to the slightly different set-ups (an underlying
point process on the compact set $K$ vs. a point process on the
noncompact space $\AA_m^d$), there are boundary effects in the
context of Theorem~\ref{thmGraphsMultivariate}, implying that the
total variation distance $\dTV(\bL_{t,a},\bM)$ is not identically
zero there.
These boundary effects are not present for $m\geq1$, which
eventually leads to the rate $O(t^{-1})$ for the
Kantorovich--Rubinstein distance
in this case.
\end{remark}

\subsection{Random polytopes with vertices on the sphere}\label{secPoissonPolytope}

\begin{figure}

\includegraphics{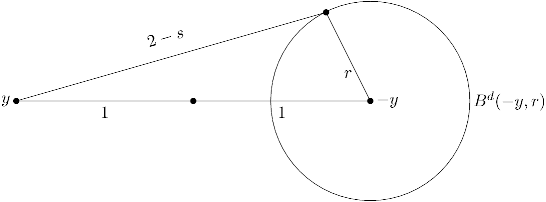}

\caption{Illustration of the argument used in the derivation of
Theorem \protect\ref{thmPoissonPolytope}.}
\label{figtriangle}
\end{figure}

Let $\SS^{d-1}$ be the unit sphere of dimension $d-1$ ($d\geq2$). Let
$\mu_t$ be a Poisson process $\eta_t$ on $\SS^{d-1}$ whose intensity
measure is a constant multiple $t\geq1$ of the normalized spherical Lebesgue
measure or a binomial process $\beta_t$ of $\lceil t\rceil$
independent and uniformly chosen points on $\SS^{d-1}$. The convex
hull $\operatorname{conv}(\mu_t)$ of $\mu_t$ is a random
polytope with vertices on $\SS^{d-1}$, and we denote by $D_t$ the
diameter of $\operatorname{conv}(\mu_t)$, that
is,
\[
D_t:=\max_{(x,y)\in\mu_{t,\neq}^2}\llVert x-y\rrVert. %
\]
More generally, define the point process of all reversed interpoint distances
by
\[
\xi_t=\frac{1}{2}\sum_{(x,y)\in\mu_{t,\neq}^2}
\delta_{2-\llVert x-y\rrVert
} . %
\]
Clearly, $D_t$ is then two minus the distance from the origin to the
closest point of~$\xi_t$. We define
\[
\bL_t(A):=\frac{1}{2} \bE\sum_{(x,y)\in\mu_{t,\neq}^2}{
\mathbf1}\bigl(t^{4/(d-1)}\bigl(2-\llVert x-y\rrVert\bigr)\in A\bigr)
, \qquad
A\subset\RR_+ \mbox{ Borel}. %
\]
Let $\chi(t):=t^2$ in the Poisson case and $\chi(t):=\lceil t\rceil
(\lceil t\rceil-1)$ in the binomial case. Applying the Mecke formula
\eqref{eqMecke} or its analogue \eqref{eqMomentMeasureBinomial} for
binomial processes, respectively, we see that
\begin{eqnarray*}
&& \bL_t\bigl([0,a]\bigr)
\\
&&\qquad =\frac{\chi(t)}{2(d\k_d)^2}\int_{\SS^{d-1}}
\int_{\SS
^{d-1}}{\mathbf1}\bigl(\llVert x-y\rrVert\geq
2-at^{-4/(d-1)}\bigr) \cH^{d-1}(\dint x) \cH^{d-1}(\dint y),
\end{eqnarray*}
where $d\k_d$ is the surface area of $\SS^{d-1}$ and $\cH^{d-1}$
stands for the $(d-1)$-dimensional Hausdorff measure.
For fixed $y\in\SS^{d-1}$, the indicator
function is one if and only if the point $x$ is contained in a certain
spherical cap $\SS^{d-1}\cap B^d(-y,r)$ centered at the antipodal
point $-y$ of $y$, whose radius $r$ has to be determined. For this, we
refer to Figure~\ref{figtriangle} and notice that $(2-s)^2+r^2=4$ so
that $r=\sqrt{4s-s^2}$. Hence the $(d-1)$-dimensional volume of $\SS
^{d-1}\cap B^d(-y,r)$ is given by
\[
(d-1)\k_{d-1}\int_0^{2s-s^2/2}
\bigl(2h-h^2\bigr)^{(d-3)/2} \,\dint h , %
\]
%

%

\noindent independently of $y$. Using the substitution
$h=2ut^{-4/(d-1)}-u^2t^{-8/(d-1)}/2$, this means
that
\begin{eqnarray*}
\bL_t\bigl([0,a]\bigr) & =&\frac{\chi(t)}{2d\k_d} (d-1)\k_{d-1}
\int_0^{2at^{-4/(d-1)}-a^2t^{-8/(d-1)}/2}\bigl(2h-h^2
\bigr)^{(d-3)/2} \,\dint h
\\
& =& \frac{\chi(t)}{2d\k_d} (d-1)\k_{d-1}
\int_0^a \bigl( 4ut^{-4/(d-1)}-u^2 t^{-8/(d-1)}
\\
&&\hspace*{97pt}{} -\bigl(2ut^{-4/(d-1)}-u^2t^{-8/(d-1)}/2
\bigr)^2 \bigr)^{(d-3)/2}
\\
&&\hspace*{93pt}{} \times\bigl(2t^{-4/(d-1)}-ut^{-8/(d-1)}\bigr)\,\dint u
\\
& =& \frac{1}{2d\k_d} \frac{\chi(t)}{t^2} (d-1)\k_{d-1}
\\
&&{}\times \int
_0^a \bigl( 4u-u^2
t^{-4/(d-1)}-t^{-4/(d-1)}\bigl(2u-u^2t^{-4/(d-1)}/2
\bigr)^2 \bigr)^{(d-3)/2}
\\
&&\hspace*{26pt}{}\times\bigl(2-ut^{-4/(d-1)}\bigr)\,\dint u .
\end{eqnarray*}
Hence we have for any Borel set $A\subset\RR_+$ that
\begin{eqnarray*}
\bL_t(A) & =& \frac{(d-1)\k_{d-1}}{2d\k_d} \frac{\chi(t)}{t^2}
\\
&&{}\times \int
_A \bigl(4u-u^2 t^{-4/(d-1)}-t^{-4/(d-1)}
\bigl(2u-u^2t^{-4/(d-1)}/2\bigr)^2
\bigr)^{(d-3)/2}
\\
&&\hspace*{24pt}{}\times\bigl(2-ut^{-4/(d-1)}\bigr)\,\dint u .
\end{eqnarray*}
The measure $\bL_t$ converges, as $t\to\infty$ and in the strong
sense, to a measure $\bM$ on $\RR_+$ given by
%
\begin{equation}
\label{eqMPoissonPolytope} \bM(A):=\frac{d-1}{d\k_d} \k
_{d-1}2^{d-3}\int
_A u^{(d-3)/2} \,\dint u ,\qquad A\subset\RR_+\mbox{ Borel} .
\end{equation}
Moreover, for any bounded Borel set $B\subset\RR_+$ there is a
constant $c_{1,B}>0$ only depending on $B$ and the space dimension $d$
such that
\[
\dTV(\bL_t\mid_B,\bM\mid_B) \leq
c_{1,B} t^{-\min\{4/(d-1),1\}}, \qquad t\geq1 . %
\]
Here, we have used that $\llvert\chi(t)/t^2-1\rrvert\leq t^{-1}$
for $t\geq1$.
Let $\widetilde{\chi}(t):=t$ in the Poisson case and $\widetilde
{\chi}(t):=\lceil t\rceil$ in the binomial case. The same arguments
as above also show that
\[
\hat r_t(B):=\sup_{x\in\SS^{d-1}}\frac{\widetilde{\chi}(t)}{d\k
_d}\int
_{\SS^{d-1}}{\mathbf1}\bigl(2-\llVert x-y\rrVert\in t^{-4/(d-1)}B
\bigr) \cH^{d-1}(\dint y)\leq c_{2,B} t^{-1}
\]
with a constant $c_{2,B}>0$ only depending on $B$ and $d$ so that $2\bL
_t(B) \hat
r_t(B)\leq2c_{2,B} \bL_t(B) t^{-1}$. Combining Corollary
\ref{corScalingLimits} and Remark~\ref{remnonUniform}(iii), we
conclude the following result. 

%
\begin{theorem}\label{thmPoissonPolytope}
Let $\zeta$ be a Poisson process on $\RR_+$ with intensity measure
given by \eqref{eqMPoissonPolytope}, and let $\xi_t$ be derived from
a Poisson process $\eta_t$ or a binomial process $\b_t$ on $\SS
^{d-1}$. Then, for any bounded Borel set $B\subset\RR_+$ there is a
constant $C_{B,d}>0$ only depending on $B$ and $d$ such that
\[
\dKR\bigl(\bigl(t^{4/(d-1)}\bullet\xi_t\bigr)
\mid_B,\zeta\mid_B \bigr)\leq C_{B,d}
t^{-\min\{4/(d-1),1\}}, \qquad t\geq1 . %
\]
In particular, for the diameter $D_t$ of the random polytope,
constructed from a Poisson process $\eta_t$ or a binomial process $\b
_t$, we have
%
\begin{eqnarray}
\bigl\llvert\bP\bigl(t^{4/(d-1)}(2-D_t)>a\bigr)-e^{-(1/(d\k_d)) \k
_{d-1}2^{d-2}a^{(d-1)/2}}
\bigr\rrvert\leq C_{a,d} t^{-\min\{4/(d-1),1\}},\nonumber
\\
\eqntext{t\geq1,}
\end{eqnarray}
with a constant $C_{a,d}>0$ only depending on $a>0$ and $d$.
\end{theorem}

%
\begin{remark}
The limiting distribution for the diameter is also derived in~\cite
{MayerMolchanov}, Theorem 5.2, and \cite{LaoMayer}, Theorem 3.1, where
the latter allows the underlying random points to have distributions
different from the uniform distribution. While the result in \cite
{MayerMolchanov} does not give any rates of convergence, in \cite
{LaoMayer}, Theorem~3.1, it has erroneously been claimed that the rate
of convergence for $D_t$ to its limiting Weibull random variable is of
order $t^{-1}$. However, in our notation the rate of convergence stated
in (2.5) in \cite{LaoMayer} concerns only the difference to a Weibull
random variable with parameter $\bL_t([0,a])$ and not to a Weibull
random variable with parameter $\bM([0,a])$ as stated by the authors.
For the difference to a Weibull random variable with parameter $\bL
_t([0,a])$, our result also yields a rate of order $t^{-1}$ since
$\mathrm{d}_{\mathrm{TV}}(\bL_t\mid_{[0,a]},\bL_t\mid_{[0,a]})=0$ in
this case.
\end{remark}


\section*{Acknowledgments}
We would like to thank an anonymous referee for valuable hints and
comments, which helped us to improve the text.

This research was initiated during the Oberwolfach mini-workshop\break 
``Stochastic Analysis for Poisson Point Processes.'' All support is
gratefully acknowledged.


%

\printaddresses
\end{document}